\documentclass[11pt]{amsart}
\usepackage{amssymb}
\usepackage{bm}
\usepackage[centertags]{amsmath}
\usepackage{amsfonts}
\usepackage{amsthm}
\usepackage{mathrsfs}
\usepackage[normalem]{ulem}
\usepackage{comment}
\usepackage{enumitem}

\usepackage{xifthen}
\usepackage{ifthenx}
\usepackage{xargs}

\usepackage[table,svgnames,x11names]{xcolor}
\usepackage[colorinlistoftodos,prependcaption,textsize=tiny,textwidth=3.5cm]{todonotes}
\presetkeys{todonotes}{fancyline}{}
\newcommandx{\change}[2][1=]{\todo[linecolor=blue,backgroundcolor=blue!25,bordercolor=blue,#1]{#2}}
\newcommandx{\changein}[2][1=]{\change[inline, caption={change}, #1]{%
    \begin{minipage}{\textwidth-20pt}#2\end{minipage}}}
\newcommandx{\todoin}[2][1=]{\todo[inline, caption={todo}, #1]{%
    \begin{minipage}{\textwidth-20pt}#2\end{minipage}}}

\newcommandx{\remove}[2][1=]{\todo[linecolor=Plum,backgroundcolor=Plum!25,bordercolor=Plum,#1]{#2}}
\newcommandx{\removein}[2][1=]{\remove[inline, caption={todo}, #1]{%
    \begin{minipage}{\textwidth-20pt}#2\end{minipage}}}

\makeatletter
\def\l@subsection{\@tocline{2}{0pt}{2.5pc}{2.5pc}{}}%
\makeatother

\usepackage{hyperref}
\hypersetup{pdftex,colorlinks=true,allcolors=black}
\usepackage{hypcap}

\DeclareFontFamily{OT1}{pzc}{}
\DeclareFontShape{OT1}{pzc}{m}{it}{<-> s * [1.10] pzcmi7t}{}
\DeclareMathAlphabet{\mathpzc}{OT1}{pzc}{m}{it}

\theoremstyle{definition}
\newtheorem{thm}{Theorem}[section]
\newtheorem{dfn}[thm]{Definition}
\newtheorem{lem}[thm]{Lemma}
\newtheorem{prop}[thm]{Proposition}

\newtheorem{rem}[thm]{Remark}

\newtheorem{ntt}[thm]{Notation}

\newtheorem*{defn*}{Definition}

\newtheorem*{thm*}{Theorem}

\newtheorem*{cor*}{Corollary}

\newtheorem*{prp*}{Proposition}

\newcommand{\N}{\mathbb{N}}

\newcommand{\trn}[1]{{\left\vert\kern-0.25ex\left\vert\kern-0.25ex\left\vert #1 
    \right\vert\kern-0.25ex\right\vert\kern-0.25ex\right\vert}}

\newcommand{\trnsmall}[1]{{\vert\kern-0.25ex\vert\kern-0.25ex\vert #1 
    \vert\kern-0.25ex\vert\kern-0.25ex\vert}}

\DeclareMathOperator{\supp}{supp}

\long\def\symbolfootnote[#1]#2{\begingroup%
\def\thefootnote{\fnsymbol{footnote}}\footnote[#1]{#2}\endgroup}

\makeatletter
\@namedef{subjclassname@2020}{%
  \textup{2020} Mathematics Subject Classification}
\makeatother

\allowdisplaybreaks

\begin{document}

\title[Orthogonal Factors of Operators on $X_{p,w}$ and $R_{\omega}^{p}$]{Orthogonal Factors of Operators on the Rosenthal $X_{p,w}$ spaces and the Bourgain-Rosenthal-Schechtman $R_{\omega}^{p}$ space}

\author[K. Konstantos]{Konstantinos Konstantos}

\address{Department of Mathematics and Statistics, York University, 4700 Keele Street, Toronto, Ontario, M3J 1P3, Canada}

\email{kostasko@yorku.ca}

\author[P. Motakis]{Pavlos Motakis}

\address{Department of Mathematics and Statistics, York University, 4700 Keele Street, Toronto, Ontario, M3J 1P3, Canada}

\email{pmotakis@yorku.ca}

\thanks{The second author was supported by NSERC Grant RGPIN-2021-03639.}




\begin{abstract} For $1<p< \infty$, we show that the Rosenthal $X_{p,w}$ spaces and the Bourgain-Rosenthal-Schechtman $R_{\omega}^{p}$ space have the factorization property and the primary factorization property.

\end{abstract}

\maketitle


\tableofcontents

\section{Introduction}

This paper is motivated by the study of conditions that when placed on a bounded linear operator $T$ on a Banach space $X$ then $T$ is a factor of the identity operator $I$ on $X$, that is there exist bounded linear operators $L,R \colon X \to X$ such that $I = LTR$. Factors of the identity have played a crucial role, e.g., in the study of decompositions of classical Banach spaces (i.e., prime and primary spaces) or the study of closed ideals in operator algebras $\mathcal{L}(X)$ (if such an ideal contains a factor of the identity then it is trivial). We place our focus on complemented subspaces of $L_{p}$, for a fixed  $1<p< \infty$. Specifically, we study operators on the Rosenthal spaces $X_{p,w}$ \cite{rosenthal:1970:Xp} and the first infinite dimensional Bourgain-Rosenthal-Schechtman  $R_{\omega}^{p}$ space \cite{bourgain:rosenthal:schechtman:1981}. Although all these spaces are isomorphic as Banach spaces (see \cite{alspach:1999} or Proposition \ref{isomorphic romegap rosenthal}), they have non-equivalent natural Schauder bases. Studying $R_{\omega}^{p}$, that has a more complicated basis, is a starting point for a similar research program for higher order  Bourgain-Rosenthal-Schechtman spaces $R_{\alpha}^{p}$, $\alpha < \omega_{1}$.

The Rosenthal $X_{p,w}$ space is the first known complemented subspace of $L_{p}$, apart from $L_{p}$, $\ell_{2}$, $\ell_{p}$, $\ell_{2} \oplus \ell_{p}$ and $\ell_{p}(\ell_{2})$. It is the closed linear span in $L_{p}$ of an independent 3-valued symmetric sequence $(f_{n})_{n=1}^{\infty}$ such that, if $p > 2$, $(\| f_{n} \|_{2} / \| f_{n} \|_{p})_{n=1}^{\infty}$ vanishes at a slow rate (when $p < 2$, we consider its conjugate exponent). The spaces $R_{\alpha}^{p}$, $\alpha < \omega_{1}$, are defined by transfinite recursion and each of them has a  Schauder basis given by a 3-valued symmetric martingale difference sequence. An appropriate cofinal subset of $\lbrace  R_{\alpha}^{p} : \alpha < \omega_{1} \rbrace$, yields uncountably many pairwise non-isomorphic complemented subspaces of $L_{p}$. We focus on $R_\omega^p$, more accurately, on its isomorphic hyperplane of all mean-zero random variables, $R_\omega^{p,0}$.  In all aforementioned cases, the complementation in $L_{p}$ is witnessed via an orthogonal projection that is also well defined on $L_{p}$.

The first result we prove is the following.
\begin{thm*} For $1<p< \infty$,  $X_{p,w}$ and $R_{\omega}^{p,0}$ with their respective standard bases have the factorization property.
\end{thm*}  
This means the following. If $X$ is one of the above spaces and $(e_{n})_{n=1}^{\infty}$ denotes its standard basis, then every bounded linear operator $T \colon X \to X$ with large diagonal, i.e., such that  $\inf_{n} \vert e_{n}^{\ast}(T(e_{n})) \vert = \delta > 0$, is a factor of the identity.

The factorization property of spaces with bases was first studied by Andrew in 1979  \cite{andrew:1979}, although it was defined slightly differently and called the perturbation property. There, he showed that $L_{p}$, $1< p < \infty$, with the Haar system has the factorization property. Several other classical spaces have been shown to have it, e.g., $H_{p}(H_{q})$, for $1 \leq p,q < \infty$ with the biparameter Haar system in \cite{laustsen:lechner:mueller:2015} and $L_{1}$ with the Haar system in \cite{LMMS:2020}.

The second result that we prove is as follows.
\begin{thm*} For $1< p < \infty$, the spaces $X_{p,w}$ and $R_{\omega}^{p,0}$ have the primary factorization property.
\end{thm*} 

The above means that for every bounded linear operator $T \colon X \to X$, where $X$ is one of these spaces, either $T$ or $I-T$ is a factor of the identity. Note that the statement of this theorem is basis-free, and the spaces  $X_{p,w}$ and $R_{\omega}^{p,0}$ are isomorphic. However, the proof relies on the basis. Specifically, we work with $R_{\omega}^{p,0}$ and its basis as a starting point for proving the factorization property of the higher complexity spaces $R_{\alpha}^{p}$, $\alpha < \omega_{1}$.

It has long been known that Banach spaces satisfying the primary factorization property together with Pe\l czy\'nski’s accordion property are primary. Thus, it has been implicitly or explicitly long studied for classical Banach spaces, often using their standard bases. Regarding the Haar system, a starting point for this can be found in the 1973 paper of Gamlen and Gaudet \cite{gamlen:gaudet:1973}, the famous Enflo-Starbird 1979 paper on $L_{1}$ \cite{enflo:starbird:1979}, and the Enflo-Maurey approach from 1975, for $1<p < \infty$ \cite{maurey:1975:2}. More recently, in 2022, in the paper  \cite{lechner:motakis:mueller:schlumprecht:2022} the biparameter Haar system of $L_{1}(L_{p})$ was used to show primariness for this space.  

In most of the above papers given a bounded linear operator $T \colon X \to X$, it is shown that there are $L,R \colon X \to X$ and a ``simpler'' operator $S \colon X \to X$ such that $LTR$ is approximately $S$. This is achieved by using as $R$ a conditional expectation (or very similar) operator
 onto a subspace $Y$ isomorphic to $X$ via an $A \colon Y \to X$ and $L = A^{-1}R$. This process is referred to as reducing $T$ to $S$. After several successive reductions, usually, $T$ is reduced to a scalar operator, i.e., an operator $\lambda I$, for some $\lambda \in \mathbb{R}$.
 
The reason why the conditional expectation has been so fruitful is that in a probability space $(\Omega, \mathcal{A}, P )$, if $(\tilde{ h}_{I})_{I \in \mathcal{D}}$ is a distributional copy of the Haar system, (for notation, see Section \ref{NC}) then $E(f) = \sum_{I \in \mathcal{D}} \vert I \vert^{-1} \langle \tilde{h}_{I}, f \rangle \tilde{ h}_{I} $ is the conditional expectation onto the subspace of $\sigma (\tilde{h}_{I})_{I \in \mathcal{D}}$-measurable functions. Something similar is true for biparameter spaces $X(Y)$ where the biparameter Haar system is a Schauder basis (see e.g., \cite{lechner:motakis:mueller:schlumprecht:2022} and    \cite{lechner:motakis:mueller:schlumprecht:2023}). The spaces  $X_{p,w}$ and $R_{\omega}^{p,0}$ also each have a basis $(f_{n})_{n=1}^{\infty}$ such that a similar approach works replacing the conditional expectation with the operator  $P(f) = \sum_{n=1}^{\infty} \|f_{n} \|_{2}^{-2}  \langle f_{n}, f \rangle f_{n}$. Although this is not a conditional expectation, it is an orthogonal projection. For our approach, we develop a language of orthogonally complemented subspaces and orthogonal factors of operators. This was also done in \cite{alspach:1999}, but at a lesser depth.      

The proof of the factorization property of $X_{p,w}$ is the least involved part of our paper and it is thus presented first. That being said, it relies on some machinery developed in \cite{LMMS:2020}, namely of strategical reproducibility of Schauder bases. This means that in a certain two-player game, one of the players has a winning strategy allowing them to reproduce the basis and its biorthogonal sequence in the dual. Using Rosenthal's analysis of block sequences of the standard basis of $X_{p,w}$, we deduce that it is strategically reproducible and thus it satisfies the factorization property. 

The study of $R_{\omega}^{p,0}$ is much more involved and it requires repeated use of probabilistic and combinatorial arguments. The space $R_{\omega}^{p,0}$ is the closed linear span of a collection $(( h_{I}^{n})_{I \in \mathcal{D}_{n-1}})_{n \in \mathbb{N}}$, where for $n \in \mathbb{N}$, $(h_{I}^{n})_{I \in \mathcal{D}_{n-1}}$ is a distributional copy of the first $(n-1)$-levels of the Haar system and the collection $\sigma ( h_{I}^{n})_{I \in \mathcal{D}_{n-1}}$,  $n \in \mathbb{N}$, is independent. With an appropriate enumeration, this basis is a martingale difference sequence, and thus, by Burkholder's inequality, it is an unconditional basis of $R_{\omega}^{p,0}$. Because each part $(h_{I}^{n})_{I \in \mathcal{D}_{n-1}}$ spans a distributional copy of $L_{p}^{n,0}$, the probabilistic techniques from  \cite{lechner:2018:1-d} apply to this setting. We refine them to construct, for every $T \colon R_{\omega}^{p,0} \to R_{\omega}^{p,0}$, an operator $D \colon R_{\omega}^{p,0} \to R_{\omega}^{p,0}$ that is diagonal with respect to the basis, and its eigenvalues are averages of the diagonal entries $(( \vert I \vert^{-1} \langle h_{I}^{n}, Th_{I}^{n} \rangle )_{I \in \mathcal{D}_{n-1}})_{n\in \mathbb{N}}$ of $T$, and show that $T$ is an approximate orthogonal factor of $D$ (or $T$ approximately orthogonally reduces to $D$). This means that, for a prechosen small $\epsilon$, there is a distributional embedding $j \colon R_{\omega}^{p,0} \to R_{\omega}^{p,0}$ such that $ \| j^{-1} P_{j( R_{\omega}^{p,0} )} T j - D \| < \epsilon$, where $P_{j( R_{\omega}^{p,0} )}$ denotes the necessarily bounded orthogonal projection onto the image of $j$. With this, the diagonal factorization property follows relatively easily from the unconditionality of the basis. The proof of the primary factorization property of  $R_{\omega}^{p,0}$ is based on an orthogonal reduction of an arbitrary diagonal operator to a scalar multiple $\lambda I$ of the identity, where $\lambda$ is an average of eigenvalues of $T$. Because approximate orthogonal factorization is transitive, every $T \colon R_{\omega}^{p,0} \to R_{\omega}^{p,0}$ is an approximate orthogonal factor of a scalar operator $\lambda I$, where $\lambda$ is in the convex hull of the diagonal entries of $T$. The definition of orthogonal factors easily then yields that $I-T$ is an approximate orthogonal factor of $(1 -\lambda)I$. The primary factorization property then follows easily. The orthogonal reduction of $D$ to $\lambda I$ is based on a probabilistic future-level stabilization argument from \cite{lechner:motakis:mueller:schlumprecht:2023} that has its roots in \cite{lechner:motakis:mueller:schlumprecht:2022}, together with some combinatorial stabilization.

The involvement of delicate probabilistic machinery is necessary for dealing with the structure of the basis $(h_{I}^{n})_{I \in \mathcal{D}_{n-1}, n \in \mathbb{N}}$. Because  $X_{p,w}$ is isomorphic to  $R_{\omega}^{p,0}$, it would have been far simpler to directly show the primary factorization property of $X_{p,w}$. However, our approach has set a path for attempting to show that all spaces $R_{\alpha}^{p}$, $\alpha < \omega_{1}$, satisfy the primary factorization property.          
 
\section{Notation and Concepts} \label{NC}

In this section, we introduce some language that will make the paper notationally easier. Moreover, we introduce the concepts of orthogonally complemented spaces (Definition \ref{defin of orthog complemented space}) and orthogonal factors (Definition \ref{defin of orthog factor}) and we give some properties related to these notions.  

By operator we mean bounded linear operator and by projection, bounded linear projection. Moreover, by subspace we mean closed linear subspace. Given $A$ a subset of a normed linear space, we denote $[A]$ the closed linear span of $A$.

Given $L_{p}(\Omega,\mathcal{A}, \mu)$, where $(\Omega,\mathcal{A}, \mu)$ is a probability space and $1 \leq p < \infty$, we denote $\mathbb{E}_{\mathcal{B}} \colon L_{p}(\Omega,\mathcal{A}, \mu) \to L_{p}(\Omega,\mathcal{A}, \mu)$ the conditional expectation operator with respect to a $\sigma$-subalgebra  $\mathcal{B}$  of $\mathcal{A}$ and if $f_{1},...,f_{n}$ are random variables on  $(\Omega,\mathcal{A}, \mu)$, we denote $\sigma(f_{1},...,f_{n})$ the $\sigma$-subalgebra of $\mathcal{A}$ generated by $f_{1},...,f_{n}$. Furthermore, if $f$ is a random variable on  $(\Omega,\mathcal{A}, \mu)$, we denote $\text{dist} (f)$ the distribution of  $f$. A space $L_{p}(\Omega,\mathcal{A}, \mu)$ is called a $L_{p}$-space. We denote $L_{p}$ the space $L_{p}[0,1]$ with the Borel $\sigma$-algebra and the Lebesque measure. 

In this paper, we consider $1< p < \infty$ and $p$ and $q$ are conjugate exponents. Moreover, denote $p^{\ast} = \max \lbrace p,q \rbrace$. For a probability space $(\Omega,\mathcal{A}, \mu)$ we consider the standard duality pairing $L_q(\Omega,\mathcal{A}, \mu)\times L_p(\Omega,\mathcal{A}, \mu)$, i.e, for $g\in L_q(\Omega,\mathcal{A}, \mu)$ and $f\in L_p(\Omega,\mathcal{A}, \mu)$ we let $\langle g,f\rangle = \int_\Omega gf d\mu$. Furthermore, if $A$ is a set, we denote $\text{card}(A)$ the cardinality of $A$ and $\chi_{A}$ the characteristic function of $A$. Finally, if $r$ is a real number, we denote $[r]$ the integer part of $r$.

{Let us recall the notion of a factor of the identity that is used in the definitions of the factorization and primary factorization property.}

\begin{dfn} Let $X$ be a Banach space and $T \colon X \to X$ be an operator.
\begin{enumerate}[label=(\alph*)]

\item We say that $T$ is a factor of the identity operator $I$ on $X$  if there are operators $A,B \colon X  \to  X$ with $I=ATB$.

\item For $K>0$ we say that $T$ is a $K$-factor of the identity operator $I$ on $X$  if there are operators $A,B \colon X \to X$ with 
\begin{align*}
\| A \| \| B \| \leq K  \;\;\; \text{and} \;\;\;   I=ATB.
\end{align*}
\end{enumerate}

\end{dfn}

The following is from \cite{laustsen:lechner:mueller:2015}. { Recall that every operator on a space with a basis can be identified with an infinite scalar matrix. A basis is said to have the factorization property if every operator whose matrix representation has diagonal entries bounded away from zero is a factor of the identity.}

\begin{dfn} Let X be a Banach space with a  basis $(e_{n})$.
\begin{enumerate}[label=(\alph*)]

\item If an operator T on X satisfies $\inf_{n} \vert e_{n}^{\ast}(T(e_{n})) \vert > 0$, then we say that T has large diagonal.

\item  We say that the basis $(e_{n})$ has the factorization property if whenever T is an operator on X with large diagonal, then $T$ is a factor of the identity operator on X.
\end{enumerate}
\end{dfn}

\begin{dfn} \label{defin of dist isom} Let $X$ and $Y$ be subspaces of $L_{p}(\Omega,\mathcal{A}, \mu)$ and $L_{p}(\Omega^{\prime},\mathcal{A}^{\prime}, \mu^{\prime})$ respectively, where $(\Omega,\mathcal{A}, \mu)$ and $(\Omega^{\prime},\mathcal{A}^{\prime}, \mu^{\prime})$  are probability spaces. We say that $X$ and $Y$  are  distributionally isomorphic if there exists a linear isomorphism $T \colon X \to Y$ with $\text{dist} (Tx) = \text{dist} (x)$, for all $x \in X$. 
\end{dfn}

We will use the following remark later, in Remark \ref{rem for bd projection} \ref{a}. After definining the notion of orthogonally complemented subspaces, this will be used to show that a distributional copy of such a space is also orthogonally complemented.

\begin{rem}  \label{remark for extension} If $X$, $Y$ and $T \colon X \to Y$ are as in Definition \ref{defin of dist isom}, then using the multivariable characteristic function we get that there exist $\sigma$-subalgebras $\mathcal{B}$ and $\mathcal{B}^{\prime}$ of $\mathcal{A}$ and $\mathcal{A}^{\prime}$  respectively, and a distributional onto isomorphism $j \colon L_{p}(\Omega, \mathcal{B}, \mu)  \to  L_{p}(\Omega^{\prime},  \mathcal{B}^{\prime}, \mu^{\prime})$ such that $j(x) = T(x)$, for all $x \in X$. 
\end{rem}

In the following definition, we introduce the notion of an orthogonally complemented space. This concept was introduced for $p > 2$ in \cite{alspach:1999}, but not used significantly.
\begin{dfn} \label{defin of orthog complemented space} Let $(\Omega,\mathcal{A}, \mu)$ be a probability space and  $C>0$.
\begin{enumerate}[label=(\alph*)]

\item An orthogonal projection $P \colon L_{2}(\Omega,\mathcal{A}, \mu) \to L_{2}(\Omega,\mathcal{A}, \mu)$ is called $p$-$C$-bounded if for every $f \in L_{2}(\Omega,\mathcal{A}, \mu) \cap  L_{p}(\Omega,\mathcal{A}, \mu)$ we have  
\begin{align*}
\| P(f) \|_{p} \leq  C \| f \|_{p} .
\end{align*}

\item  A subspace $X$ of $L_{p}(\Omega,\mathcal{A}, \mu)$ is called $C$-orthogonally complemented if there exists a $p$-$C$-bounded projection $P$ such that $X$ is the closure of $P ( L_{2}(\Omega,\mathcal{A}, \mu) \cap  L_{p}(\Omega,\mathcal{A}, \mu))$ with respect to $\| \cdot \|_{p}$. In this case, we denote $P_{X}$ the projection on $L_{p}(\Omega,\mathcal{A}, \mu)$ onto $X$, which is given by extending $P \restriction_{L_{2}(\Omega,\mathcal{A}, \mu) \cap  L_{p}(\Omega,\mathcal{A}, \mu)}$ to $L_{p}(\Omega,\mathcal{A}, \mu)$. Note $\| P_{X} \| \leq C$.
\end{enumerate}

\end{dfn}

{Next we recall the definition of a martingale difference sequence. The basis of the Bourgain-Rosenthal-Schechtman  $R_{\omega}^{p}$ space is a martingale difference sequence.}

\begin{dfn} Let $(f_{n})$ be a sequence in $L_{p}(\Omega, \mathcal{A}, \mu)$, where $(\Omega,\mathcal{A}, \mu)$ is a probability space. The sequence  $(f_{n})$ is called a martingale difference sequence if  $\mathbb{E}(f_{n} | f_{1},...,f_{n-1})=0$ , for every $n  \geq 2$.
\end{dfn}

Let us state some properties of orthogonally complemented spaces used in the rest of the paper.  
\begin{rem} \label{rem for bd projection}  Let $(\Omega,\mathcal{A}, \mu)$ be a probability space and  $C>0$.
\begin{enumerate}[label=(\alph*)]

\item \label{a} Every distributional copy of a $C$-orthogonally complemented space is $C$-orthogonally complemented. This follows from Remark  \ref{remark for extension} combined with   a conditional expectation operator.

\item \label{b} By duality and self-adjointness, if $P$ is a  $p$-$C$-bounded projection, then  $P$ is a  $q$-$C$-bounded projection.

\item If $X$ is a subspace of $L_{p}(\Omega,\mathcal{A}, \mu)$ and it is  $C$-orthogonally complemented, then the dual space $X^{\ast}$ of $X$ is $C$-isomorphic to the closure of $P ( L_{2}(\Omega,\mathcal{A}, \mu) \cap  L_{q}(\Omega,\mathcal{A}, \mu))$ with respect to $\| \cdot \|_{q}$, which is $C$-orthogonally complemented, from \ref{b}.

\item \label{d} If $(f_{n})$ is a $\lbrace -1,0,1 \rbrace $-valued symmetric martingale difference sequence and $X = [(f_{n})]$ in  $L_{p}(\Omega,\mathcal{A}, \mu)$ is $C$-orthogonally complemented, then 
\begin{align*}
P_{X}(f) = \sum_{n=1}^{\infty} \langle \| f_{n} \|_{q}^{-1} f_{n}, f \rangle   \| f_{n} \|_{p}^{-1} f_{n}  .
\end{align*}
This is because $\| f_{n} \|_{p} \| f_{n} \|_{q} = \| f_{n} \|_{2}^{2} $ (see \cite[page 280]{rosenthal:1970:Xp}). The convergence of the above series is in $\| \cdot \|_{p}$, by the Burkholder inequality.

\end{enumerate}

\end{rem}

In the following definition, we introduce the notion of an approximate orthogonal factor. This is a special case of an operator being $T$ an approximate factor of an operator $S$, and the factoring operators relate to an orthogonal projection.

\begin{dfn} \label{defin of orthog factor} Let $C>0$ and $\epsilon >0$. Let $X$ and $Y$ be subspaces of $L_{p}(\Omega,\mathcal{A}, \mu)$ and $L_{p}(\Omega^{\prime},\mathcal{A}^{\prime}, \mu^{\prime})$ respectively, where $(\Omega,\mathcal{A}, \mu)$ and $(\Omega^{\prime},\mathcal{A}^{\prime}, \mu^{\prime})$  are probability spaces, and assume that $Y$ is $C$-orthogonally  complemented. Moreover, let $T \colon X \to X$ and $S \colon Y \to Y$ be operators. We say that $T$ is an orthogonal factor of $S$ with error $\epsilon$ if there exists a distributional embedding $j \colon Y \to X$ such that $\| j^{-1} P_{j(Y)} T j -S \| < \epsilon$. When the error is $\epsilon =0$, we will simply say that $T$ is an orthogonal factor of $S$.

\end{dfn}

\begin{rem} Let $T$ and $S$ be as in Definition \ref{defin of orthog factor}. If $T$ is an orthogonal factor of $S$ with error $\epsilon >0$, then $I_{X} - T$ is an orthogonal factor of $I_{Y} - S$ with error $\epsilon >0$, where  $I_{X}$ and $I_{Y}$ are the identity operators on $X$ and $Y$ respectively. This is because $P_{j(Y)} j = j$ and therefore 
\begin{align*}
\| j^{-1} P_{j(Y)} (I_{X} - T) j - (I_{Y} - S) \| < \epsilon.
\end{align*}

\end{rem}

The following explains the relation between orthogonal factors and projectional factors. The notion of projectional factor can be found in \cite[Definition 2.1]{lechner:motakis:mueller:schlumprecht:2022}, but we do not explicitly use it here.

\begin{rem} Let $C>0$ and $\epsilon >0$. Let $X$ be a subspace of $L_{p}(\Omega,\mathcal{A}, \mu)$, where $(\Omega,\mathcal{A}, \mu)$ is a probability space, and assume that $X$ is $C$-orthogonally complemented. Moreover, let $T,S \colon X \to X$ be operators. If $T$ is an orthogonal factor of $S$ with error $\epsilon >0$, then $T$ is a $C$-projectional  factor of $S$ with error $\epsilon >0$. Note $C$ is the orthogonal complementation constant of $X$.

\end{rem}

In a certain sense, being an approximate orthogonal factor is a transitive property. 

\begin{prop} \label{transitivity}

Let $T \colon X \to X$, $S \colon Y \to Y$ and $R \colon Z \to Z$ be operators, where $X$, $Y$ and $Z$ are subspaces of possibly different $L_{p}$-spaces, and assume that $Y$ is $C$-orthogonally complemented and $Z$ is $D$-orthogonally complemented. Furthermore, assume that $T$ is an orthogonal factor of $S$ with error $\epsilon_{1}$ and $S$ is an orthogonal factor of $R$ with error $\epsilon_{2}$. Then $T$ is an orthogonal factor of $R$ with error $ D \epsilon_{1} + \epsilon_{2}$.
\begin{proof}
By the assumptions we have that 
\begin{align*}
\| j_{1}^{-1} P_{j_{1}(Y)} T j_{1}  -S \| < \epsilon_{1}
\end{align*} 
and
\begin{align*}
\| j_{2}^{-1} P_{j_{2}(Z)} S j_{2}  - R \| < \epsilon_{2},
\end{align*} 
where $j_{1} \colon Y \to X$ and $j_{2} \colon Z \to Y$ are distributional embeddings.
Put $j=j_{1}j_{2}$ and observe that  
\begin{align*}
j_{2}^{-1} P_{j_{2}(Z)} j_{1}^{-1} P_{j_{1}(Y)}T j_{1} j_{2} = j^{-1} P_{j(Z)} T j .
\end{align*}
Hence, 
\begin{align*}
& \| j^{-1} P_{j(Z)} T j - R \| \leq \| j^{-1} P_{j(Z)} T j - j_{2}^{-1} P_{j_{2}(Z)} S j_{2} \| + \| j_{2}^{-1} P_{j_{2}(Z)} S j_{2}  - R \| \\
& \leq  \| j_{2}^{-1} \|  \| P_{j_{2}(Z)} \|  \| j_{1}^{-1} P_{j_{1}(Y)} T j_{1} - S \| \| j_{2} \| +  \| j_{2}^{-1} P_{j_{2}(Z)} S j_{2}  - R \| \\
& \leq D \epsilon_{1} + \epsilon_{2}.
\end{align*}
Therefore, $T$ is an orthogonal factor of $R$ with error $ D \epsilon_{1} + \epsilon_{2}$.
\end{proof}

\end{prop}

\section{The factorization property of the basis of the $X_{p,w}$ space}

In this section, we recall the definition of the $X_{p,w}$ space (Definition \ref{Xpw space}) and the Rosenthal $X_{p}$ space (Definition \ref{Xp space}). In addition, we state some results related to these spaces that can be found in \cite{rosenthal:1970:Xp}, and we deduce that the Rosenthal $X_{p}$ space is isomorphic to a $7.35(p/\log(p))$-orthogonally complemented subspace of $L_{p}$ (Theorem \ref{Xp is orth complem}). Finally, using the notion of a strategically reproducible basis that can be found in \cite{LMMS:2020}, we prove that the unit vector basis of the  $X_{p,w}$ space has the factorization property (Theorem \ref{big theorem for Xp}). Although all $X_{p,w}$ spaces, for $w$ satisfying the $(\ast)$ property, are mutually isomorphic, their bases are not equivalent. But, we show that all of them have the factorization property. Note $X_{p}$ does not have a standard basis.

\begin{dfn} \label{Xpw space} Let $2<p< \infty$ and $w=(w_{n})$ be a sequence of positive reals satisfying the $(\ast)$ property : 
\begin{align*}
w_{n} \rightarrow 0 \;\;\; \text{and} \;\;\; \sum_{n=1}^{\infty}  w_{n}^{2p/(p-2)} = \infty.
\end{align*}
Define 
\begin{align*}
X_{p,w} = \lbrace (x_{n}) \subset \mathbb{R} : (x_{n}) \in \ell_{p}  \; \text{and} \; (x_{n} w_{n}) \in \ell_{2} \rbrace
\end{align*}
and 
\begin{align*}
\| (x_{n}) \|_{X_{p,w}} = \max \lbrace \| (x_{n}) \|_{p}, 
 \| (x_{n}w_{n}) \|_{2} \rbrace.
\end{align*} 
The $(X_{p,w}, \| \cdot \|_{X_{p,w}})$ space is a Banach space with unconditional and shrinking basis the unit vector basis $(e_{n})$.
\end{dfn}

The following can be found in \cite[Theorem 13]{rosenthal:1970:Xp}.

\begin{thm}[] Let $2<p< \infty$ and $w, w^{\prime}$ be sequences of positive reals satisfying the $(\ast)$ property. Then $X_{p,w}$ is isomorphic to $X_{p,w^{\prime}}$.
\end{thm}

\begin{dfn} \label{Xp space} Let $2<p< \infty$. The Rosenthal $X_{p}$ space is an $X_{p,w}$ space, where $w=(w_{n})$ is a sequence of positive reals satisfying the $(\ast)$ property. 
\end{dfn}

Recall that, for brevity, we denote $L_p[0,1]$ by $L_p$. It is proved in \cite[Corollary, page 282]{rosenthal:1970:Xp} that for $2 < p < \infty$, the Rosenthal $X_{p}$ space is isomorphic to a complemented subspace of $L_{p}$. Furthermore, we have the following:

\begin{thm} \label{Xp is orth complem} Let $2<p< \infty$. The Rosenthal $X_{p}$ space is isomorphic to a $7.35(p/\log(p))$-orthogonally complemented subspace of $L_{p}$.
\end{thm}

This is not stated in \cite{rosenthal:1970:Xp}, but it can be easily deduced as follows: The Rosenthal $X_{p}$ space is isomorphic to the closed linear span of a $3$-valued symmetric independent sequence $(f_{n})$ in $L_{p}$ and this closed linear span is $7.35(p/\log(p))$-orthogonally complemented subspace of $L_{p}$ because the restriction of the orthogonal projection $P \colon L_{2} \to L_{2}$, defined by
\begin{align*}
P(f) = \sum_{n=1}^{\infty} \langle f_{n},f \rangle \| f_{n} \|_{2}^{-2} f_{n}, 
\end{align*}
to $L_{p}$ is a projection onto the closed linear span of $(f_{n})$ in $L_{p}$ of norm at most $K_{p}$ (\cite[Theorem 4]{rosenthal:1970:Xp}). The constant $K_{p}$ is at most $7.35(p/\log(p))$  from \cite[Theorem 4.1]{JSZ:1985}.

The following can be founded in \cite[Lemma 7]{rosenthal:1970:Xp}.

\begin{lem} \label{blocks} Let $2<p< \infty$ and $w$ be a sequence of positive reals satisfying the $(\ast)$ property.  Let $(E_{j})$ be a sequence of disjoint finite subsets of the positive integers. For each j, put
\begin{align*}
f_{j}  =  \sum_{n \in E_{j}} w_{n}^{2/(p-2)}e_{n},
\end{align*}
\begin{align*}
\beta_{j}  =  \left( \sum_{n \in E_{j}} w_{n}^{2p /(p-2)} \right)^{(p-2)/2p} 
\end{align*}
and
\begin{align*}
\tilde{{f}_{j}} & =   \| f_{j} \|_{p}^{-1} f_{j}.
\end{align*}
Let $Y$ denote the closed linear span of $(\tilde{{f}_{j}})$ in   $X_{p,w}$. Then  $(\tilde{{f}_{j}})$ is an unconditional basis for $Y$, isometrically equivalent to the unit vector basis of $X_{p,(\beta_{j})}$ and $Y$ is a 1-complemented subspace of  $X_{p,w}$ with projection defined by 
\begin{align*}
P (x) = \sum_{j=1}^{\infty} \langle \| f_{j} \|_{p} \| f_{j} \|^{-2}_{2} f_{j},x \rangle  \tilde{{f}_{j}}.
\end{align*}

\end{lem}

The following two definitions are from  \cite{LMMS:2020}. See Remark \ref{2-player game} below for an intuitive explanation of them.

\begin{dfn} Let $(x_{n})$ and $(y_{n})$ be Schauder basic sequences in (possibly different) Banach spaces and $C \geq 1$. We say that $(x_{n})$ and $(y_{n})$ are impartially $C$-equivalent if for any finite choice of scalars $(a_{n}) \in  c_{00}$ we have 
\begin{align*}
\frac{1}{\sqrt{C}} \| \sum_{n=1}^{\infty} a_{n} y_{n} \| \leq \| \sum_{n=1}^{\infty} a_{n} x_{n} \| \leq  \sqrt{C}\| \sum_{n=1}^{\infty} a_{n} y_{n} \| . 
\end{align*}

\end{dfn}

\begin{dfn} \label{def of sr} Assume that X is a Banach space with a basis $(e_{n})$ which is  unconditional and shrinking. Let $(e_{n}^{\ast}) \subset X^{\ast}$ be the corresponding coordinate functionals. We say that $(e_{n})$  is strategically reproducible if the following condition is satisfied for some $C \geq 1$:

\begin{align*}
\forall n_{1} \in \mathbb{N} & \exists b_{1} \in \langle  e_{n} : n \geq n_{1} \rangle, & \exists b^{\ast}_{1} \in \langle  e^{\ast}_{n} : n \geq n_{1} \rangle\\
\forall n_{2} \in \mathbb{N} & \exists b_{2} \in \langle  e_{n} : n \geq n_{2} \rangle, & \exists b^{\ast}_{2} \in \langle  e^{\ast}_{n} : n \geq n_{2} \rangle\\
\forall n_{3} \in \mathbb{N} & \exists b_{3} \in \langle  e_{n} : n \geq n_{3} \rangle, & \exists b^{\ast}_{3} \in \langle  e^{\ast}_{n} : n \geq n_{3} \rangle\\
\span\vdots
\end{align*}
such that 

\begin{enumerate}[label=(\alph*)] 

\item $(b_{k}) \text{ is impartially C-equivalent to} (e_{k}),$
\item $(b^{\ast}_{k}) \text{ is impartially C-equivalent to} (e^{\ast}_{k}),$ and 
\item $b^{\ast}_{k}(b_{l}) = \delta_{k,l}  \forall k,l \in \mathbb{N}$
\end{enumerate}
\end{dfn} 

\begin{rem} \label{2-player game} The condition in Definition \ref{def of sr} can be interpreted  as player (II) having a winning strategy in a two-player game: \\
We fix $C \geq 1$. Player (I) chooses $n_{1} \in \mathbb{N}$, then player (II) chooses 
\begin{align*}
b_{1} \in \langle  e_{n} : n \geq n_{1} \rangle \;\;\; \text{and} \;\;\;  b^{\ast}_{1} \in \langle  e^{\ast}_{n} : n \geq n_{1} \rangle.
\end{align*}
They repeat the moves infinitely many times, obtaining, for every $k \in \mathbb{N}$, numbers $n_{k}$, and vectors $b_{k}$ and $b^{\ast}_{k}$. Player (II) wins if they were able to choose the
sequences $(b_{k}) \subset X$ and $(b^{\ast}_{k}) \subset X^{\ast}$ such that (a), (b) and (c) are satisfied. Thus, the basis $(e_{n})$ is strategically reproducible if and only if for some $C \geq 1$ player (II) has a winning strategy.
 
\end{rem}

The following theorem will be used inside the proof of the main theorem of this section (Theorem \ref{big theorem for Xp}) and it can be found in \cite[Theorem 3.12]{LMMS:2020}.

\begin{thm} \label{sr} Let $X$ be a  Banach space with an unconditional and shrinking basis $(e_{n})$. If the basis $(e_{n})$ is strategically reproducible, then the basis has the factorization property.
\end{thm}

\begin{thm} \label{big theorem for Xp} Let $2<p< \infty$. For every $w=(w_{n})$ satisfying the $(\ast)$ property, the unit vector basis of the $X_{p,w}$ space has the factorization property.
\begin{proof} By Theorem \ref{sr} it is enough to show that the basis $(e_{n})$ for the $X_{p,w}$ space is strategically reproducible. In the two-player game of Remark \ref{2-player game} we assume the role of player (II). 

Let $\epsilon >0$. In each round, after player (I) chooses $n_{k} \in \mathbb{N}$, we choose a finite subset $E_{k}$ of the positive integers such that $\min(E_{k}) > n_{k}$, $\max(E_{k-1}) < \min(E_{k}) $, where $E_{0} = \emptyset$, and such that putting 
\begin{align*}
\beta_{k}  = \left( \sum_{n \in E_{k}} w_{n}^{2p/(p-2)} \right)^{(p-2)/2p} 
\end{align*}
we have 
\begin{align*}
 w_{k} \leq \beta_{k} \leq \sqrt{1+ \epsilon} w_{k} . 
\end{align*} 
Now, put 
\begin{align*}
b_{k} &= \sum_{n \in E_{k}} w_{n}^{2/(p-2)}e_{n}, 
\end{align*}
\begin{align*}
\tilde{b}_{k} &= \| b_{k} \|_{p}^{-1} b_{k} 
\end{align*}
and
\begin{align*}
b^{\ast}_{k} &= \sum_{n \in E_{k}} \| b_{n} \|_{p} \| b_{n} \|^{-2}_{2} w^{2/(p-2)}_{n} e^{\ast}_{n}. 
\end{align*}
We show that this choice satisfies Definition \ref{def of sr}. We observe that $\tilde{b}_{k} \in \langle  e_{n} : n \geq n_{k} \rangle$ and $b^{\ast}_{k} \in \langle  e^{\ast}_{n} : n \geq n_{k} \rangle$.
Moreover, put $\beta = (\beta_{k})$. Hence $(\tilde{b}_{k})$ is isometrically equivalent to the usual basis for the  $X_{p,\beta}$ space by Lemma \ref{blocks}. Now, we will show that the usual basis for $X_{p, \beta}$  is impartially $(1+ \epsilon)$-equivalent to the usual basis for $X_{p,w}$. 

Indeed, let $a_{1},...,a_{n} \in \mathbb{R}$. We have
\begin{align*}
\|  \sum_{k=1}^{n} a_{k}e_{k} \|_{X_{p, \beta}} = \max \lbrace \|  \sum_{k=1}^{n} a_{k}e_{k} \|_{p}, \|  \sum_{k=1}^{n} a_{k} \beta_{k} e_{k} \|_{2} \rbrace 
\end{align*}
and
\begin{align*}
\sum_{k=1}^{n} a_{k}e_{k} \|_{X_{p, w}} = \max \lbrace \|  \sum_{k=1}^{n} a_{k}e_{k} \|_{p}, \|  \sum_{k=1}^{n} a_{k} w_{k} e_{k} \|_{2} \rbrace.
\end{align*}
But, since for every $k \in \mathbb{N}$,  $w_{k} \leq \beta_{k} \leq \sqrt{1+ \epsilon} w_{k}$, we get  
\begin{align*}
\|  \sum_{k=1}^{n} a_{k} \beta_{k} e_{k} \|_{2} \leq  \sqrt{1 + \epsilon} \|  \sum_{k=1}^{n} a_{k} w_{k} e_{k} \|_{2} 
\end{align*} 
and
\begin{align*}    
\|  \sum_{k=1}^{n} a_{k} w_{k} e_{k} \|_{2} \leq  \|  \sum_{k=1}^{n} a_{k} \beta_{k} e_{k} \|_{2}. 
\end{align*}
Hence, 
\begin{align*}
\|  \sum_{k=1}^{n} a_{k}e_{k} \|_{X_{p, w}} \leq \|  \sum_{k=1}^{n} a_{k}e_{k} \|_{X_{p, \beta}}  \leq  \sqrt{1+ \epsilon}  \|  \sum_{k=1}^{n} a_{k}e_{k} \|_{X_{p, w}},
\end{align*}
as desired. Thus, $(\tilde{{b}}_{k})$ is impartially $(1+ \epsilon)$-equivalent to the usual basis for $X_{p,w}$. Let $Y$ be the closed linear span of the  sequence $(\tilde{{b}}_{k})$  in  $X_{p, w}$. By Lemma \ref{blocks} we know that $Y$ is 1-complemented subspace of $X_{p, w}$ with   
\begin{align*}
P (x) = \sum_{k=1}^{\infty} \langle \| b_{k} \|_{p} \| b_{k} \|^{-2}_{2} b_{k},x \rangle  \tilde{{b}}_{k} 
\end{align*}
as the norm-1 projection onto $Y$.  Since  $(\tilde{{b}}_{k})$ is impartially $(1+ \epsilon)$-equivalent to the usual basis for $X_{p,w}$, there exists an isomorphism $T \colon Y \to  X_{p,w} $  with  $T(\tilde{{b}}_{k}) = e_{k}$ and $\| T \| \leq \sqrt{1+ \epsilon}$. We observe that  $b^{\ast}_{k} = e^{\ast}_{k} T P$ and  $b^{\ast}_{k}( \tilde{{b}}_{l} ) = \delta_{k,l}$. Now, it remains to prove that $( b^{\ast}_{k}) $ is impartially $(1+ \epsilon)$-equivalent to $(e^{\ast}_{k})$. 

Indeed, let $a_{1},...,a_{n} \in \mathbb{R}$. We observe that 
\begin{align*}
\| \sum_{k=1}^{n} a_{k} b^{\ast}_{k} \|  = \| \sum_{k=1}^{n} a_{k} e^{\ast}_{k} T P \|=  \| (TP)^{\ast} \left( \sum_{k=1}^{n} a_{k} e^{\ast}_{k} \right) \| \leq  \sqrt{1+ \epsilon} \|  \sum_{k=1}^{n} a_{k} e^{\ast}_{k} \|. 
\end{align*}
For the other inequality, we have
\begin{align*}
\|  \sum_{k=1}^{n} a_{k} e^{\ast}_{k} \| =  \left( \sum_{k=1}^{n} a_{k} e^{\ast}_{k} \right) (x) , 
\end{align*}
where $x$ belongs to the unit sphere of $ X_{p,w}$ and $x = \sum_{k=1}^{\infty} \lambda_{k}e_{k}.$ We have 
\begin{align*}
\| \sum_{k=1}^{n} a_{k} e^{\ast}_{k} \| =  \sum_{k=1}^{n} a_{k} \lambda_{k} = \left(    \sum_{k=1}^{n} a_{k} b^{\ast}_{k} \right) \left( \sum_{k=1}^{n} \lambda_{k} \tilde{{b}_{k}}  \right) \leq  \sqrt{1+ \epsilon} \| \sum_{k=1}^{n} a_{k} b^{\ast}_{k} \|, 
\end{align*}
as desired.
\end{proof}
\end{thm}

It is easy to see that a reflexive Banach space with Schauder basis $(e_{n})$ has the factorization property if and only if the dual space $X^{\ast}$ with the biorthogonals $(e_{n}^{\ast})$ has the factorization property. For $1<q<2$, Rosenthal (see \cite[page 280]{rosenthal:1970:Xp}) defines $X_{q,w}$ to be the dual space $X_{p,w}^{\ast}$  of  $X_{p,w}$. Therefore, we have the following.
\begin{thm} For $1<q<2$ and $w=(w_{n})$ satisfying the $(\ast)$ property with respect to $p$, the Rosenthal $X_{q,w}$ space with the biorthogonals of the  unit vector basis of the $X_{p,w}$ space has the factorization property.

\end{thm}

From the paragraph below Theorem \ref{Xp is orth complem}, it follows that, for $1<q<2$, the  Rosenthal $X_{q,w}$ space is isomorphic to the closed linear span in $L_{q}$ of an independent  sequence of 3-valued symmetric random variables $(f_{n})$.

\section{The Bourgain-Rosenthal-Schechtman  $R_{\alpha}^{p}$ spaces } 

In this section we recall the Bourgain-Rosenthal-Schechtman  $R_{\alpha}^{p}$ spaces (Definition \ref{defin of BRS Spaces}) using the disjoint and independent sums (Definition \ref{sum}), and we prove that the $R_\alpha^p$ spaces are orthogonally complemented (Theorem \ref{Orth Compl spaces the BRS}) using facts from \cite{bourgain:rosenthal:schechtman:1981}. We describe a specific distributional copy of $R_{\omega}^{p}$ using the finite dimensional $L_{p}^{n}$ spaces (Remark \ref{import rem for R omega space}) and establish an explicit formula for the orthogonal projection onto this copy (Remark \ref{important for boundedness of E}) expressed by the martingale different sequence  basis of this copy. We will work with this fixed copy of $R_\omega^p$ for the rest of the paper.

\subsection{Haar system notation} We recall basic notation for the standard Haar system and finite-dimensional $L_p$-spaces.

\begin{ntt} A dyadic interval of $[0,1]$ is one of the form $[\frac{i-1}{2^{n}},  \frac{i}{2^{n}})$, where $n \in \mathbb{N} \cup \lbrace 0 \rbrace$ and $1 \leq  i  \leq 2^{n}$. If $I= [\frac{i-1}{2^{n}},  \frac{i}{2^{n}})$, then denote $I^{+} = [ \frac{2(i-1)}{2^{n+1}},  \frac{2i -1}{2^{n+1}})$, $I^{-} = [ \frac{2i-1}{2^{n+1}},  \frac{2i}{2^{n+1}})$ 
and $\pi(I^{+}) = \pi(I^{-}) = I$. Denote the collection of the dyadic intervals by $\mathcal{D}$ and denote 
\begin{align*}
\mathcal{D}^{n} & = \lbrace I \in \mathcal{D} : \vert I \vert = 2^{-n} \rbrace ,
\end{align*}
and
\begin{align*}
\mathcal{D}_{n} & = \bigcup_{k=0}^{n} \mathcal{D}^{k}. 
\end{align*}

Define the linear order on $\mathcal{D}^{n}$ as follows: If $I,J \in \mathcal{D}^{n}$ with $I = [\frac{i-1}{2^{n}},  \frac{i}{2^{n}})$ and $J = [\frac{j-1}{2^{n}},  \frac{j}{2^{n}})$, where $1 \leq  i,j  \leq 2^{n}$, then we say that $I<J$ if and only if $i<j$. Furthermore, define the linear order on $\mathcal{D}_{n}$ as follows: If $I,J \in \mathcal{D}_{n}$ with $I \in \mathcal{D}^{k}$ and $J \in \mathcal{D}^{m}$, where $0 \leq k, m \leq n$ then we say that $I<J$ if and only if $k<m$ or if $k = m$ then $I<J$ in the linear order of $\mathcal{D}^{k}$.

\end{ntt}
{In the following notation we recall the definition of the standard Haar system and the notation  of  the finite dimensional $L_{p}$-spaces that we will use to describe the Bourgain-Rosenthal-Schechtman $R_{\omega}^{p}$ space.} 
\begin{ntt}\label{Haar system} Recall that $L_{p} = L_{p}[0,1]$. We recall the definition of the Haar system $(h_{I})_{ I \in \mathcal{D}}$ on $L_{p}$. For $I \in \mathcal{D}$ we denote $h_{I}$ by
\begin{align*}
h_{I} = \chi_{I^{+}} - \chi_{I^{-}}.
\end{align*}
By the Paley-Burkholder inequality the Haar system $(\chi_{[0,1]}) \cup(h_{I})_{ I \in \mathcal{D}}$ is an  unconditional basis of $L_{p}$ for $1 < p < \infty$ and therefore $L_{p} = [(\chi_{[0,1]}) \cup  (h_{I})_{ I \in \mathcal{D}}]$. Denote
\begin{align*}
L_{p}^{0} & = \lbrace f \in L_{p} : \mathbb{E}(f) = 0 \rbrace. 
\end{align*}
We observe that 
\begin{align*}
L_{p}^{0} = [(h_{I})_{ I \in \mathcal{D}}].
\end{align*} 
Let $n \in \mathbb{N}$. Denote
\begin{align*}
L_{p}^{n} = \langle (\chi_{[0,1]}) \cup (h_{I})_{I \in \mathcal{D}_{n -1}}  \rangle = \langle \chi_{I} : I \in \mathcal{D}^{n} \rangle.
\end{align*} 
and
\begin{align*}
L_{p}^{n,0} = \lbrace f \in L_{p}^{n} : \mathbb{E}(f) = 0 \rbrace .
\end{align*}
Similarly
\begin{align*}
L_{p}^{n,0} = \langle (h_{I})_{I \in \mathcal{D}_{n-1}}  \rangle.
\end{align*} 
\end{ntt}

\subsection{The $R_\alpha^p$ spaces} 

In this subsection, we recall the definition of the Bourgain-Rosenthal-Schechtman  $R_{\alpha}^{p}$ spaces (Definition \ref{defin of BRS Spaces}) and we prove that the $R_\alpha^p$ spaces are orthogonally complemented (Theorem \ref{Orth Compl spaces the BRS}). 

In the following definition, we recall two types of sums used in the construction of the $R_{\alpha}^{p}$  (Definition \ref{defin of BRS Spaces}), i.e., disjoint sums and independent sums.
\begin{dfn} \phantom{A} \label{sum} 

\begin{enumerate}[label=(\alph*)]

\item Let $B$  be a subspace of $L_{p}(\Omega, \mathcal{A}, \mu)$. Define the disjoint sum $ (B  \oplus B)_{p}$ of $B$ to be the  subspace
\[ 
\begin{split}
\Big\{  b(\omega, \epsilon)\in L_p(\Omega\times\{0,1\}) : &\text{ there exists }    b_{\epsilon} \in B \text{ with }  b(\omega, \epsilon) =   b_{\epsilon} (\omega)\\
&\text{for all } \omega \in \Omega, \;\epsilon = 0 \text{ or } 1 \Big\} ,
\end{split}
\]
of $L_p(\Omega\times\{0,1\})$, where $\lbrace 0,1 \rbrace$ is endowed with the fair probability measure assigning mass $\frac{1}{2}$ to each 0 and 1.

\item Let $B_{n}$ be a subspace of $L_{p}(\Omega_{n}, \mathcal{A}_{n}, \mu_{n})$, where $n=1,2,3,...$. Define the independent sum of the spaces $B_{n}$ as follows: Take the product probability space $\prod_{n} \Omega_{n}$ and  for each $n$, denote $\overline{B}_{n}$ the subspace 
\[ \begin{split} 
\Bigg\{ b \in L_{p} \left( \prod_{n} \Omega_{n}  \right) : \; \text{there exists} \; f \in B_{n} \; \text{with} \; b(\omega) = f(\omega_{n}) \\ 
\text{for all} \; \omega \in \prod_{n} \Omega_{n} \Bigg\},
\end{split}
\]
of $L_{p} \left(  \prod_{n} \Omega_{n}  \right)$, that is  $\overline{B}_{n}$ is a distributional copy of $B_{n}$ depending only on the $n$-th coordinate. Then the independent sum  $( \sum_{n} B_{n} ) _{\text{Ind},p}$  of the spaces $B_{n}$ is $ [\bigcup_{n} \overline{B}_{n}]$ in $L_{p} \left(  \prod_{n} \Omega_{n}  \right)$.

\end{enumerate}

\end{dfn}

\begin{rem} \label{rem for ind sum} Let $B_{n}$ be a subspace of $L_{p}(\Omega_{n}, \mathcal{A}_{n}, \mu_{n})$, where $n=1,2,3,...$ and let $L_{p}(\Omega,\mathcal{A}, \mu)$. Assume that we can find independent $\sigma$-subalgebras $\mathcal{S}_{n}$ of $\mathcal{A}$ and distributional copies $X_{n}$ of $B_{n}$ such that if $f \in X_{n}$, then $f$ is $\mathcal{S}_{n}$-measurable, that is the spaces  $X_{n}$ are distributional independent copies of $B_{n}$ in $L_{p}(\Omega,\mathcal{A}, \mu)$, then the independent sum  $( \sum_{n} B_{n} ) _{\text{Ind},p}$ is distributionally isomorphic to the space $[\bigcup_{n} X_{n}]$ in  $L_{p}(\Omega,\mathcal{A}, \mu)$.
\end{rem}

The Bourgain-Rosenthal-Schechtman  $R_{\alpha}^{p}$ spaces are defined recursively, starting with the spaces of constant functions. Disjoint sums are used in the successor case and  independent sums are used in the limit case.

\begin{dfn} \label{defin of BRS Spaces} Define $R_{0}^{p}$ to  be the space of the  constant functions of $L_{p}\lbrace 0,1 \rbrace$ and let $\beta$ be an ordinal with $0< \beta < \omega_{1}$. Suppose that we have defined the $R_{a}^{p}$ space  for all $\alpha < \beta$.
\begin{enumerate}[label=(\alph*)]

\item If $\beta = \alpha + 1$, then $R_{\beta}^{p}  = ( R_{a}^{p} \oplus R_{a}^{p})_{p}$.

\item If $\beta$ is a limit ordinal, then $  R_{\beta}^{p} =  ( \sum_{\alpha < \beta} R_{\alpha}^{p} ) _{\text{Ind},p}$

\end{enumerate}

\end{dfn}
Below, we track some details from \cite{bourgain:rosenthal:schechtman:1981} that establish the orthogonal complementation of the Bourgain-Rosenthal-Schechtman  $R_{\alpha}^{p}$ spaces.

\begin{thm}\label{Orth Compl spaces the BRS} Let $\alpha < \omega_{1}$. The $R_{\alpha}^{p}$ space is $(p^{\ast} -1)^{2} (\frac{p^{\ast}}{2})^{3/2}$-orthogonally complemented.
\end{thm}

\begin{proof} The proof is based on technical details from \cite{bourgain:rosenthal:schechtman:1981}. We give a description that avoids some of them, but the unfamiliar reader is invited to skip this proof. For $1<p<\infty$, Bourgain, Rosenthal, and Schechtman define an explicit auxiliary subspace $X_\mathcal{D}^p$ of $L_p\{0,1\}^\mathcal{D}$, for some countable set $\mathcal{D}$. The specific form of $\mathcal{D}$ and $X_\mathcal{D}^p$ is not relevant to this proof, but, for context, we mention it is the closed linear span, with respect to $\|\cdot\|_p$, of a vector space $X$ of measurable functions that is the same for all values of $p$. For more details on $X_\mathcal{D}^p$, the interested reader may refer to \cite[Section 3]{bourgain:rosenthal:schechtman:1981}. According to \cite[Theorem 3.1]{bourgain:rosenthal:schechtman:1981}, $X_\mathcal{D}^p$ is complemented in $L_p\{0,1\}^\mathcal{D}$. More precisely, it is stated in \cite[page 222, paragraph above equation 3.8]{bourgain:rosenthal:schechtman:1981} that the projection $P_0$ onto $X_\mathcal{D}^p$ is in fact orthogonal. A precise formula for $P_0$ is given, and it is independent of $p$. The estimate $\|P_0\| \leq K_p^2(p/2)^{3/2}$ (where $K_p = (p-1)$ is the constant in Burkholder's inequality, see, e.g. \cite[Theorem 8.6]{AdamOse:2012}) is given in \cite[bottom of page 222]{bourgain:rosenthal:schechtman:1981} only for the case $p>2$. However, the ``selfadjointness'' of $P_0$ yields that the same estimate is true for $1<p<2$ replacing $p$ with $p^*$, i.e., $\|P_0\| \leq (K_{p^{\ast}})^2(p^{\ast}/2)^{3/2}$, for $1<p<\infty$.

Fixing $\alpha<\omega_1$, it is shown in \cite[Lemma 3.9]{bourgain:rosenthal:schechtman:1981} that there exists a distributional copy $X_{T_\alpha}^p$ of $R_\alpha^p$ in $X_\mathcal{D}^p$. In \cite[line 1 of proof of Lemma 3.6]{bourgain:rosenthal:schechtman:1981} it is stated that $X_{T_\alpha}^p$ is complemented in $X_\mathcal{D}^p$ via a conditional expectation operator $\mathbb{E}_\alpha$. Therefore, $\mathbb{E}_\alpha P_0$ is a $p$-$C$-bounded orthogonal projection on $L_p\{0,1\}^\mathcal{D}$ onto $X_{T_\alpha}^p$, for $C =  (K_{p^{\ast}})^2(p^{\ast}/2)^{3/2}$. By Remark \ref{rem for bd projection} \ref{a}, $R_{\alpha}^{p}$ is $(p^{\ast} -1)^{2} (\frac{p^{\ast}}{2})^{3/2}$-orthogonally complemented.
\end{proof}

\subsection{ The $R_\omega^p$ space and its basis }

{In this subsection we describe a specific distributional copy of $R_{\omega}^{p}$ using the finite dimensional  $L_{p}^{n}$ spaces and establish an explicit formula for the orthogonal projection onto this copy.}

\begin{rem} \label{import rem for R omega space} By Definition \ref{defin of BRS Spaces} we have that $R_{\omega}^{p} =  ( \sum_{n=0} R_{n}^{p} ) _{\text{Ind},p}$. Because, $R_0^p = \langle\chi_{[0,1]}\rangle$ it is easy to see by induction that for every $n \in \mathbb{N}$ the $R_{n}^{p}$ space is a distributional copy of $L_{p}^{n}$. Hence, $R_{\omega}^{p}$ is a distributional copy of $( \langle\chi_{[0,1]}\rangle\oplus \sum_{n=1} L_{p}^{n} ) _{\text{Ind},p}$. We may therefore identify $R_\omega^p$ with this copy, that is,
\[R_\omega^p = ( \langle\chi_{[0,1]}\rangle\oplus \sum_{n=1} L_{p}^{n} ) _{\text{Ind},p}.\]

\end{rem}

\begin{ntt} Denote
\begin{align*}
R_{\omega}^{p,0} & = \lbrace f \in R_{\omega}^{p} : \mathbb{E}(f) = 0 \rbrace .
\end{align*}
\end{ntt}

\begin{rem}
\label{important rem for R omega space 2}
It is more convenient for the purpose of this paper to work with $R_{\omega}^{p,0}$ instead of $R_{\omega}^{p}$. Note that
\[( \langle \chi_{[0,1]} \rangle\oplus \sum_{n=1} L_{p}^{n} ) _{\text{Ind},p} = ( \sum_{n=1} L_{p}^{n,0} ) _{\text{Ind},p} + \langle \chi_{[0,1]} \rangle.\]
So, by Remark \ref{import rem for R omega space},
\[R_{\omega}^{p,0} = ( \sum_{n=1} L_{p}^{n,0} ) _{\text{Ind},p}.\]
\end{rem}

The following is proved in \cite{alspach:1999} (see the penultimate paragraph of page 11 until the first seven lines of page 12).

\begin{prop}
\label{isomorphic romegap rosenthal}
The Bourgain-Rosenthal-Schechtman space $R_\omega^{p,0}$ and the Rosenthal space $X_p$ are isomorphic.
\end{prop}

Note that $R_\omega^p$ is isomorphic to its hyperplanes as an infinite-dimensional complemented subspace of $L_p$. In particular, $R_\omega^p\simeq R_\omega^{p,0}\simeq X_p$.

We will use the following linearly ordered set to enumerate the basis of $R_{\omega}^{p}$.  
\begin{ntt} Denote 
\begin{align*}
\mathcal{D}_{\omega} & = \lbrace (n,I) : n \in \mathbb{N}, I \in \mathcal{D}_{n-1} \rbrace
\end{align*} 
and define the linear order on $\mathcal{D}_{\omega}$ as follows: $(m,J) \prec (n,I)$ if and only if $m<n$ or if $m=n$, then $J < I$ in the linear order of  $\mathcal{D}_{m-1}$. {The set $\mathcal{D}_{\omega}$ can be visualized as disjoint copies of the sets $\mathcal{D}_{n-1}$, $n\in\N$.}
\end{ntt}
The subscript $\omega$ for the set $\mathcal{D}_{\omega}$ indicates the space $R_{\omega}^{p}$. Similar sets can be used for the $R_{\alpha}^{p}$ spaces.
{We now define the basis of $R_{\omega}^{p}$ and in Remark \ref{MDS the sequence} we observe that it is a martingale difference sequence.}
\begin{ntt} \label{Importan notation for the basis of R}
We fix a specific basis for $R_{\omega}^{p}$ inside $L_p$ as follows. For every $n$ we fix  $(h_{I}^{n})_{I \in \mathcal{D}_{n-1}} $ a distributional copy of the basis $(h_{I})_{I \in \mathcal{D}_{n-1}}$ of $L_{p}^{n,0}$ in $L_{p}$ such that the  $\sigma$-subalgebras $\sigma(  \lbrace  h_{I}^{n} : I \in \mathcal{D}_{n-1} \rbrace )$, $n \in \mathbb{N}$, are independent. Then, by Remark \ref{import rem for R omega space},  $R_{\omega}^{p}$ is distributionally isomorphic to $[ (\chi_{[0,1]}) \cup (h_{I}^{n})_{(n,I) \in \mathcal{D}_{\omega}} ]$ in  $L_{p}$ and, by Remark \ref{important rem for R omega space 2}, $R_{\omega}^{p,0}$ is distributionally isomorphic to $[(h_{I}^{n})_{(n,I) \in \mathcal{D}_{\omega}} ]$ in  $L_{p}$. The sequence $ (h_{I}^{n})_{ (n,I) \in \mathcal{D}_{\omega} }$ consists of  $\lbrace -1, 0, 1 \rbrace $-valued symmetric random variables. For the rest of the paper, we will identify  $R_\omega^p$ and $R_\omega^{p,0}$ with their aforementioned distributional copies in $L_p$, i.e.,
\begin{align*}
R_{\omega}^{p} = [ (\chi_{[0,1]}) \cup (h_{I}^{n})_{(n,I) \in \mathcal{D}_{\omega}} ]
\end{align*}
and
\begin{align*}
R_{\omega}^{p,0} = [ (h_{I}^{n})_{(n,I) \in \mathcal{D}_{\omega}} ].
\end{align*}
\end{ntt}

\begin{rem}\label{MDS the sequence}  The sequence $(h_{I}^{n})_{(n,I) \in \mathcal{D}_{\omega}}$ described in Notation \ref{Importan notation for the basis of R} is a martingale difference sequence with respect to the linear order of $\mathcal{D}_{\omega}$. {This follows from the fact that  $(h_{I}^{n})_{I \in \mathcal{D}_{n-1}} $, $n \in \mathbb{N}$, are independent and for every $n \in \mathbb{N}$  we have that $(h_{I}^{n})_{I \in \mathcal{D}_{n-1}} $ is a distributional copy of the martingale difference sequence basis $(h_{I})_{I \in \mathcal{D}_{n-1}}$ of $L_{p}^{n,0}$}. Therefore, $(h_{I}^{n})_{(n,I) \in \mathcal{D}_{\omega}}$ is an unconditional basis of $R_{\omega}^{p,0}$ by the Burkholder inequality.
\end{rem}

{In the following remark we explain that $R_{\omega}^{p,0}$ is orthogonally complemented.}

\begin{rem} \label{Ort Compl for Proj} By Theorem \ref{Orth Compl spaces the BRS} the $R_{\omega}^{p}$ space is $(p^{\ast} -1)^{2} (\frac{p^{\ast}}{2})^{3/2}$-orthogonally complemented {in $L_p$ via some projection $P$. Because the space $L_p^0$ is  2-orthogonally completed via the projection $Qf =  f-\mathbb{E}(f)$ and $R_\omega^{p,0}$ is a subspace of $L_p^0$, it follows that $PQ$ is a $p$-$C$-bounded orthogonal projection onto $R_\omega^{p,0}$, for $C = 2(p^{\ast} -1)^{2} (\frac{p^{\ast}}{2})^{3/2}$. Therefore, $R_\omega^{p,0}$ is $2(p^{\ast} -1)^{2} (\frac{p^{\ast}}{2})^{3/2}$-orthogonally complemented.}
\end{rem}

{We now give an explicit formula for the projection onto $R_{\omega}^{p,0}$.}
\begin{rem} \label{important for boundedness of E}

By Remarks \ref{MDS the sequence}, \ref{Ort Compl for Proj} and \ref{rem for bd projection} \ref{d} we have that the projection $P_{R_{\omega}^{p,0}} \colon  L_{p}  \to   L_{p}$ onto  $R_{\omega}^{p,0}$ is   
\begin{align*}
P_{R_{\omega}^{p,0}} (f) = \sum_{n=1}^{\infty} \sum_{I \in \mathcal{D}_{n-1}} \vert I \vert^{-1} \langle h_{I}^{n}, f \rangle h_{I}^{n}.
\end{align*}
By Remarks \ref{rem for bd projection} \ref{a} and \ref{d}, if $(b_I^n)_{(n,I)\in\mathcal D_\omega}$ is a distributional copy of $(h_I^n)_{(n,I)\in\mathcal D_\omega}$, then replacing $h_I^n$ with $b_I^n$ in the above formula yields the orthogonal projection onto $[(b_I^n)_{(n,I)\in\mathcal D_\omega}]$.
\end{rem}

\section{Orthogonal Reduction to a diagonal operator}

In this section, we approximately orthogonally reduce an arbitrary operator on $R_{\omega}^{p,0}$ to a diagonal operator on $R_{\omega}^{p,0}$ (Theorem \ref{thm1}) using Lemma \ref{lem 4}. This is based on probabilistic arguments first used by Lechner in \cite{lechner:2018:1-d}. Finally, using Theorem \ref{thm1} we get that an arbitrary operator on $R_{\omega}^{p,0}$ with large diagonal is a factor of the identity operator on $R_{\omega}^{p,0}$ (Theorem \ref{large diagonal theorem}).

\begin{ntt} \label{ntt 2} The following will be used inside the proof of Theorem \ref{thm1}. We give this notation now to discuss the proof.

\begin{enumerate}[label=(\alph*)] 

\item Let $(b^{n}_{I})_{(n,I) \in \mathcal{D}_{\omega}}$ be a distributional copy of $(h^{n}_{I})_{(n,I) \in \mathcal{D}_{\omega}}$  in $R_{\omega}^{p,0}$ and let $Y$ be the closed linear span of $(b^{n}_{I})_{(n,I) \in \mathcal{D}_{\omega}}$. Define the distributional isomorphism $j \colon R_{\omega}^{p,0}  \to Y$ with 
\begin{align*}
j(h_{I}^{n}) = b^{n}_{I}
\end{align*}
and the projection  $E \colon R_{\omega}^{p,0}  \to R_{\omega}^{p,0}$ onto $Y$ with 
\begin{align*}
E(f) = \sum_{n=1}^{\infty} \sum_{I \in \mathcal{D}_{n-1}} \vert I \vert^{-1} \langle b_{I}^{n}, f \rangle b_{I}^{n} ,
\end{align*}
which is bounded by Remark \ref{important for boundedness of E}.

\item Let $T \colon R_{\omega}^{p,0}  \to R_{\omega}^{p,0}$ be an operator and denote  
\begin{align*}
\mathcal{D}(T)= \lbrace \vert I \vert^{-1} \langle h_{I}^{n}, Th_{I}^{n} \rangle :   (n,I) \in \mathcal{D}_{\omega} \rbrace 
\end{align*}
the diagonal elements of T with respect to the basis  $(h^{n}_{I})_{(n,I) \in \mathcal{D}_{\omega}}$ of   $R_{\omega}^{p,0}$. Define $\mathcal{A}(\mathcal{D}(T))$ to be the collection of all averages of finite subsets of $\mathcal{D}(T)$.

\end{enumerate}

\end{ntt}

We describe now in broad strokes the general idea of the proof of Theorem \ref{thm1}. Let $T \colon R_{\omega}^{p,0}  \to R_{\omega}^{p,0}$ be an arbitrary operator. Using an induction on the linear order of $\mathcal{D}_\omega$, for every $n \in  \mathbb{N}$ we construct blocks $(b^{n}_{I})_{I \in \mathcal{D}_{n-1}}$ of the basis $(h^{n}_{I})_{I \in \mathcal{D}_{n-1}}$ of $L_{p}^{n,0}$ in some $L_{p}^{N(n),0}$, $N(n)> n$, with $b_{I}^{n} = \sum _{K \in \mathcal{B}_{I}^{n}} \theta_{K} h_{K}^{N(n)}$, where $( \theta_{K})_{K \in \mathcal{B}_{I}^{n}} \subset \lbrace -1, 1 \rbrace$, such that 
\begin{align*}
&\supp(b_{[0,1]}^{n}) = [0,1] \\
&\supp( b_{I^{+}}^{n}) = [ b_{I}^{n} = 1] \\
&\supp( b_{I^{-}}^{n}) = [ b_{I}^{n} = -1],
\end{align*} 
for every $n \in  \mathbb{N}$. From this we obtain that $(b^{n}_{I})_{(n,I) \in \mathcal{D}_{\omega}}$ is a distributional copy of $(h^{n}_{I})_{(n,I) \in \mathcal{D}_{\omega}}$. We construct the block sequence $(b^{n}_{I})_{(n,I) \in \mathcal{D}_{\omega}}$ such that 
\begin{align*}
\vert \langle b_{J}^{m}, Tb_{I}^{n} \rangle \vert < \epsilon_{J,I}^{m,n},
\end{align*} 
for every $(m,J) \neq (n,I) \in \mathcal{D}_{\omega}$, and 
\begin{align*}
\vert \langle b_{I}^{n}, Tb_{I}^{n} \rangle - \sum_{ K \in \mathcal{B}_{I}^{n}} \langle h_{K}^{N(n)}, T h_{K}^{N(n)} \rangle \vert < \epsilon_{I}^{n},
\end{align*}
for every $(n,I) \in \mathcal{D}_{\omega}$, where $\epsilon_{J,I}^{m,n}$ and $\epsilon_{I}^{n}$ are small errors. From this we obtain that $T$ is an approximate orthogonal factor of a diagonal operator $R$, that is 
$$\| j^{-1}ETj -R \| < \epsilon,$$ 
where $j$ and $E$ are as in Notation \ref{ntt 2} and the diagonal operator $R$ is defined by  
$$R(h_{I}^{n}) = \left( \vert I \vert ^{-1} \sum_{K \in \mathcal{B}_{I}^{n}} \langle h_{K}^{N(n)}, Th_{K}^{N(n)} \rangle \right) h_{I}^{n}$$ 
and has diagonal element in $\mathcal{A}(\mathcal{D}(T))$.

In each inductive step, the choice of the block vector $b^{n}_{I}$ is performed probabilistically, and for this purpose, we give the following notation.
 
\begin{ntt} \label{Random Variables} Let $T \colon R_{\omega}^{p,0}  \to R_{\omega}^{p,0}$ be an operator, $f \in L_{q}^{0}$ and $x \in R_{\omega}^{p,0}$. Moreover, let $M,N \in \N$ with $M \geq N$ and $\mathcal{B} \subseteq \mathcal{D}^{N}$. Consider the probability space $\{-1,1\}^\mathcal{B}$ with the fair probability measure. Define the random vector
\begin{align*}
b_{\mathcal{B},M}(\theta) = \sum_{K \in \mathcal{B}} \theta_{K} h_{K}^{M},
\end{align*}
where $\theta = (\theta_{K})_{K \in \mathcal{B}}  \in \lbrace -1, 1 \rbrace ^{\mathcal{B}}$ and  the following random variables
\begin{align*}
 Y_{\mathcal{B},M,f}(\theta) &= \langle f,  b_{\mathcal{B},M}(\theta) \rangle, 
\end{align*}
\begin{align*}
 W_{\mathcal{B},M,x}(\theta) &= \langle  b_{\mathcal{B},M}(\theta), x \rangle 
\end{align*}
and
\begin{align*}
 Z_{\mathcal{B},M,T}(\theta) &= \langle b_{\mathcal{B},M}(\theta), Tb_{\mathcal{B},M}(\theta) \rangle - \sum_{K \in \mathcal{B}} \langle h_{K}^{M}, Th_{K}^{M} \rangle.
\end{align*}
The vector $b_{\mathcal{B},M}$ is in $L_q$ and $L_p$ and thus can occupy either position in $\langle\cdot,\cdot\rangle$.
\end{ntt}

In the above notation, the random vector $b_{\mathcal{B},M}$ will be used in the inductive step of the previously described proof to choose an appropriate vector $b^{n}_{I}$. The random variables $Y_{\mathcal{B},M,f}$, $W_{\mathcal{B},M,x}$, and $Z_{\mathcal{B},M,T}$, represent the interaction of this vector with previously constructed vectors via $T$, and they are used to achieve the small errors $\epsilon_{J,I}^{m,n}$ and $\epsilon_{I}^{n}$. Next, we find estimates for the variance of these random variables. We will use these estimates in the proof of Lemma \ref{lem 4}.
\begin{lem} \label{lem1} Let $M,N \in \N$ with $M \geq N$, $\mathcal{B} \subseteq \mathcal{D}^{N}$ and $f \in L_{q}^{0}$. We have  
\begin{align*}
\mathbb{E}(Y_{\mathcal{B},M,f}) = 0
\end{align*}
and
\begin{align*}
\mathrm{Var}(Y_{\mathcal{B},M,f}) \leq \| f \|_{q}^{2}  \vert \cup \mathcal{B} \vert^{1/p} 2^{-N/p}.
\end{align*}
\begin{proof} We have 
\begin{align*}
Y_{\mathcal{B},M,f}(\theta) = \langle f,  b_{\mathcal{B},M}(\theta) \rangle = \sum_{K \in \mathcal{B}} \theta_{K} \langle f,  h_{K}^{M} \rangle.
\end{align*} 
So, since
\begin{align*}
\mathbb{E}(Y_{\mathcal{B},M,f}) = \sum_{K \in \mathcal{B}} \mathbb{E} (\theta_{K}) \langle f,  h_{K}^{M} \rangle =0,
\end{align*}
we get that  
\begin{align*}
\mathrm{Var}(Y_{\mathcal{B},M,f}) = \mathbb{E}(Y^{2}_{\mathcal{B},M,f}).
\end{align*}
But,
\begin{align*}
 Y^{2}_{\mathcal{B},M,f} (\theta) &= (\sum_{K \in \mathcal{ B}} \theta_{K} \langle f,  h_{K}^{M} \rangle) (\sum_{L \in \mathcal{B}} \theta_{L} \langle f,  h_{L}^{M} \rangle)  \\
&= \sum_{K,L \in \mathcal{B}} \theta_{K} \theta_{L} \langle f,  h_{K}^{M} \rangle \langle f,  h_{L}^{M} \rangle.
\end{align*}
Now, taking the expectation we get 
\begin{align*}
\mathbb{E}(Y^{2}_{\mathcal{B},M,f}) = \sum_{K,L \in \mathcal{B}} \mathbb{E}(\theta_{K} \theta_{L}) \langle f,  h_{K}^{M} \rangle \langle f,  h_{L}^{M} \rangle .
\end{align*}
However, $ \mathbb{E}(\theta_{K} \theta_{L}) =1$ if and only if $K=L$ and, otherwise, $ \mathbb{E}(\theta_{K} \theta_{L}) =0$ . Therefore,
$$\mathrm{Var}(Y_{\mathcal{B},M,f}) = \sum_{K \in B}  \langle f,  h_{K}^{M} \rangle^{2}.$$ 
Set $a_{f,K} = \langle f, h_{K}^{M} \rangle$ and  we observe that $\vert  a_{f,K} \vert \leq \| f \|_{q}  \vert K \vert ^{1/p}$. Hence, 
\begin{align*}
&\mathrm{Var}(Y_{\mathcal{B},M,f}) = \langle f, \sum_{K \in \mathcal{B}} a_{f,K} h_{K}^{M} \rangle \leq  \| f \|_{q} \| \sum_{K \in \mathcal{B}} a_{f,K} h_{K}^{M}  \|_{p} \\
&\leq  \| f \|_{q} \max_{K \in \mathcal{B}} \vert  a_{f,K} \vert  \| \sum_{K \in \mathcal{B}} h_{K}^{M}  \|_{p} = \| f \|_{q}  \max_{K \in \mathcal{B}} \vert  a_{f,K} \vert  \vert \cup \mathcal{B} \vert ^{1/p}  \\
&\leq  \| f \|_{q}^{2}   \vert \cup \mathcal{B} \vert ^{1/p} \max_{K \in \mathcal{B}} \vert  K \vert ^{1/p}  = \| f \|_{q}^{2}  \vert \cup \mathcal{B} \vert^{1/p}  2^{-N/p}.
\end{align*}
\end{proof}
\end{lem}

\begin{lem} \label{lem2}  Let $M,N \in \N$ with $M \geq N$, $\mathcal{B} \subseteq \mathcal{D}^{N}$ and $x \in R_{\omega}^{p,0}$. We have
\begin{align*}
\mathbb{E}(W_{\mathcal{B},M,x}) = 0
\end{align*}
and
\begin{align*}
\mathrm{Var}(W_{\mathcal{B},M,x}) \leq \| x \|_{p}^{2}  \vert \cup \mathcal{B} \vert^{1/q} 2^{-N/q}.
\end{align*}
\begin{proof} The proof is almost identical to that of Lemma \ref{lem1}.
\end{proof}
\end{lem}

\begin{lem} \label{lem 3} Let $T \colon R_{\omega}^{p,0}  \to  R_{\omega}^{p,0}$ be an operator, $M,N \in \N$ with $M \geq N$ and $\mathcal{B} \subseteq \mathcal{D}^{N}$. We have 
\begin{align*}
\mathbb{E}(Z_{\mathcal{B},M,T}) = 0
\end{align*}
and
\begin{align*}
\mathrm{Var}(Z_{\mathcal{B},M,T}) \leq 2 \| T \|^{2} \vert \cup \mathcal{B} \vert ^{(1/p) +1} 2^{-N/q}.
\end{align*}
\begin{proof} We have 
\begin{align*}
Z_{\mathcal{B},M,T}(\theta) &= \langle b_{\mathcal{B},M}(\theta), Tb_{\mathcal{B},M}(\theta) \rangle - \sum_{K \in \mathcal{B}} \langle h_{K}^{M}, Th_{K}^{M} \rangle \\
&= \sum _{K \neq L \in \mathcal{B}} \theta_{K} \theta_{L} \langle h_{K}^{M}, Th_{L}^{M} \rangle.
\end{align*}
So, 
\begin{align*}
\mathbb{E}(Z_{\mathcal{B},M,T}) = \sum _{K \neq L \in \mathcal{B}} \mathbb{E}(\theta_{K} \theta_{L}) \langle h_{K}^{M}, Th_{L}^{M} \rangle = 0 
\end{align*}
because $K \neq L$. Thus 
\begin{align*}
\mathrm{Var}(Z_{\mathcal{B},M,T}) = \mathbb{E}(Z^{2}_{\mathcal{B},M,T}). 
\end{align*}
But,
\begin{align*}
Z^{2}_{\mathcal{B},M,T}(\theta) &= \left( \sum _{K_{1} \neq L_{1} \in \mathcal{B}} \theta_{K_{1}} \theta_{L_{1}} \langle h_{K_{1}}^{M}, Th_{L_{1}}^{M} \rangle \right)  \left( \sum _{K_{2} \neq L_{2} \in \mathcal{B}} \theta_{K_{2}} \theta_{L_{2}} \langle h_{K_{2}}^{M}, Th_{L_{2}}^{M} \rangle \right) \\ 
&=\sum _{
 K_{1} \neq L_{1} \in \mathcal{B} ,
K_{2} \neq L_{2} \in \mathcal{B} } 
\theta_{K_{1}} \theta_{L_{1}}\theta_{K_{2}}\theta_{L_{2}} \langle h_{K_{1}}^{M}, Th_{L_{1}}^{M} \rangle \langle h_{K_{2}}^{M}, Th_{L_{2}}^{M} \rangle.
\end{align*}
Now, taking the expectation we get 
\begin{align*}
\mathbb{E}(Z^{2}_{\mathcal{B},M,T}) = \sum _{ K_{1} \neq L_{1} \in \mathcal{B} ,
K_{2} \neq L_{2} \in \mathcal{B} } 
\mathbb{E} ( \theta_{K_{1}} \theta_{L_{1}}\theta_{K_{2}}\theta_{L_{2}}) \langle h_{K_{1}}^{M}, Th_{L_{1}}^{M} \rangle \langle h_{K_{2}}^{M}, Th_{L_{2}}^{M} \rangle.
\end{align*}
However, $\mathbb{E} ( \theta_{K_{1}} \theta_{L_{1}}\theta_{K_{2}}\theta_{L_{2}}) =1 $ if and only if $K_{1} = K_{2}, L_{1} = L_{2}$ and $K_{1} = L_{2}, K_{2} = L_{1}$, otherwise, $\mathbb{E} ( \theta_{K_{1}} \theta_{L_{1}}\theta_{K_{2}}\theta_{L_{2}}) = 0 .$ Therefore, 
\begin{align*}
\mathrm{Var}(Z_{\mathcal{B},M,T}) =  \sum _{
 K \neq L \in \mathcal{B} }  \langle h_{K}^{M}, Th_{L}^{M} \rangle  \langle h_{L}^{M}, Th_{K}^{M} \rangle  +  \sum _{
 K \neq L \in B }  \langle h_{K}^{M}, Th_{L}^{M} \rangle^{2} .
\end{align*}
Set 
\begin{align*}
A =  \sum _{
 K \neq L \in \mathcal{B} }  \langle h_{K}^{M}, Th_{L}^{M} \rangle  \langle h_{L}^{M}, Th_{K}^{M} \rangle
\end{align*} 
and
\begin{align*}
B &=  \sum _{
 K \neq L \in \mathcal{B} }  \langle h_{K}^{M}, Th_{L}^{M} \rangle^{2}.
\end{align*}
We will find  estimations for A and B.  Indeed, set $a_{L,K} =  \langle h_{L}^{M}, Th_{K}^{M} \rangle$. We have $\vert a_{L,K} \vert \leq \| T \| \vert L \vert^{1/q} \vert K \vert^{1/p}.$ Hence, 
\begin{align*}
A &= \sum_{K \in \mathcal{B}} \langle h_{K}^{M}, \sum_{L  \neq K \in \mathcal{B}} a_{L,K} Th_{L}^{M} \rangle \\
& \leq \| T \| \sum_{K \in \mathcal{B}} \vert K \vert^{1/q} \| \sum_{L  \neq K \in \mathcal{B}} a_{L,K} h_{L}^{M}  \|_{p} \\
& \leq \| T \| \sum_{K \in \mathcal{B}} \vert K \vert^{1/q} \max_{L \neq K \in \mathcal{B}} \vert a_{L,K} \vert   \| \sum_{L  \neq K \in \mathcal{B}} h_{L}^{M}  \|_{p} \\
& \leq  \| T \| \sum_{K \in \mathcal{B}} \vert K \vert^{1/q} \max_{L \neq K \in \mathcal{B}} \vert a_{L,K} \vert   \| \sum_{L  \in \mathcal{B}} h_{L}^{M}  \|_{p} \\ 
& \leq  \| T \|^{2} \vert \cup \mathcal{B} \vert^{1/p} \sum_{K \in \mathcal{B}} \vert K \vert \max_{L \neq K \in \mathcal{B}} \vert L \vert^{1/q} \\
 &= \| T \|^{2} \vert \cup \mathcal{B} \vert^{(1/p) +1} 2^{-N/q}.
\end{align*}
The proof of the estimation for B is almost identical to that of A and we have  $B \leq \| T \|^{2} \vert \cup \mathcal{B} \vert^{(1/p) +1} 2^{-N/q}$. Therefore, 
\begin{align*}
\mathrm{Var}(Z_{\mathcal{B},M,T}) \leq 2 \| T \|^{2} \vert \cup \mathcal{B} \vert ^{(1/p) +1}  2^{-N/q}.
\end{align*}
\end{proof}
\end{lem}
Now, we use the previous probabilistic estimates to determine a value for a random vector. This will be used in the inductive step of the proof of Theorem \ref{thm1}.
\begin{lem} \label{lem 4} Let $T \colon R_{\omega}^{p,0}  \to R_{\omega}^{p,0}$ be an operator. Let $(n,I) \in \mathcal{D}_{\omega}$, $M,N \in \mathbb{N}$ with $M \geq N$ and $\mathcal{B} \subseteq \mathcal{D}^{N}$ with $\vert \cup \mathcal{B} \vert = \vert I \vert$. Let $x_{(m,J)} \in R_{\omega}^{p,0}$ and $f_{(m,J)} \in L_{q}^{0}$, where $(m,J) \prec (n,I)$ such that 
\begin{align*}
\| x_{(m,J)} \| & \leq \vert J \vert^{1/p} 
\end{align*}
and
\begin{align*}
 \| f_{(m,J)} \| & \leq \vert J \vert^{1/q}.
\end{align*}
For $\theta = (\theta_{K})_{K \in \mathcal{B}}  \in \lbrace -1, 1 \rbrace ^{\mathcal{B}}$ we  define the random variables 
\begin{align*}
Z(\theta) = Z_{\mathcal{B},M,T} (\theta) =  \langle b_{\mathcal{B},M}(\theta), Tb_{\mathcal{B},M}(\theta) \rangle - \sum_{K \in \mathcal{B}} \langle h_{K}^{M}, Th_{K}^{M} \rangle ,
\end{align*}
\begin{align*}
Y_{(m,J)}(\theta) = Y_{\mathcal{B},M,T^{\ast}f_{(m,J)}}(\theta)  = \langle T^{\ast}f_{(m,J)},   b_{\mathcal{B},M}(\theta) \rangle ,  \;\;\;   (m,J) \prec (n,I)  
\end{align*}
and
\begin{align*}
W_{(m,J)}(\theta) = W_{\mathcal{B},M,Tx_{(m,J)}}(\theta) = \langle  b_{\mathcal{B},M}(\theta), Tx_{(m,J)}  \rangle , \;\;\; (m,J) \prec (n,I).
\end{align*}
Let $\eta >0$ and $\eta_{(m,J)}>0$, where $(m,J) \prec (n,I) $. If 
\begin{equation} \label{star}
N >   p^{\ast}  \left( 2 \log \| T \| + \log  \left( \frac{2^{2n+3}}{\eta^{2}} + \sum_{(m,J) \prec (n,I)} \frac{1}{\eta^{2}_{(m,J)}} \right) \right) 
\end{equation}
 
then there exists $\theta^{0} = (\theta_{K}^{0})_{K \in \mathcal{B}}$ such that 
\begin{align*}
\vert Z(\theta^{0}) \vert < \eta,
\end{align*}
\begin{align*}
\vert Y_{(m,J)} (\theta^{0}) \vert < \eta, \;\;\; (m,J) \prec (n,I)
\end{align*}
and
\begin{align*} 
\vert W_{(m,J)} (\theta^{0}) \vert < \eta_{(m,J)}, \;\;\; (m,J) \prec (n,I).
\end{align*}
\begin{proof} To prove the existence of $\theta^{0}$, we define the following events:
\begin{align*}
\mathcal{Z} = \lbrace \theta : \vert Z(\theta) \vert > \eta \rbrace ,
\end{align*}
\begin{align*}
\mathcal{Y}(m,J) = \lbrace \theta : \vert Y_{(m,J)} (\theta) \vert > \eta \rbrace, \qquad (m,J) \prec (n,I)
\end{align*}
and 
\begin{align*}
\mathcal{W}(m,J) = \lbrace \theta : \vert W_{(m,J)} (\theta) \vert > \eta_{(m,J)} \rbrace, \qquad (m,J) \prec (n,I).
\end{align*}
Let
\[
\mathcal{A} = \mathcal{Z} \cup\Big(\bigcup_{(m,J) \prec (n,I)} \mathcal{Y}(m,J)\Big)\cup\Big(  \bigcup_{(m,J) \prec (n,I)} \mathcal{W}(m,J)\Big).
\]

Taking the probability of the above events, we observe that 
\begin{align*}
\mathbb{P}( \mathcal{A})\leq  \mathbb{P}(\mathcal{Z}) + \sum_{(m,J) \prec (n,I)} \mathbb{P} ( \mathcal{Y}(m,J)) + \sum_{(m,J) \prec (n,I)} \mathbb{P} ( \mathcal{W}(m,J) )
\end{align*}
\begin{align*}
\leq  \frac{1}{\eta^{2}} \mathrm{Var}(Z) + \frac{1}{\eta^{2}} \sum_{(m,J) \prec (n,I)} \mathrm{Var} (Y_{(m,J)}) + \sum_{(m,J) \prec (n,I)} \frac{1}{\eta^{2}_{(m,J)}} \mathrm{Var} (W_{(m,J)}),
\end{align*}
by the Chebyshev inequality. Now, by Lemmas \ref{lem1}, \ref{lem2} and \ref{lem 3}, we get 
\begin{align*}
&\mathbb{P}( \mathcal{A}) \leq   \frac{2}{\eta^{2}} \| T \|^{2} \vert I \vert^{(1/p)+1} 2^{-N/q}
 +  \frac{1}{\eta^{2}}  \sum_{(m,J) \prec (n,I)} \| T \|^{2} \vert J \vert^{2/q} \vert I \vert ^{1/p} 2^{-N/p} \\ 
&+ \sum_{(m,J) \prec (n,I)} \frac{1}{\eta^{2}_{(m,J)}} \| T \|^{2} \vert J \vert^{2/p} \vert I \vert ^{1/q} 2^{-N/q}  \\
& \leq \left( \frac{2^{2n+3}}{\eta^{2}} + \sum_{(m,J) \prec (n,I)} \frac{1}{\eta^{2}_{(m,J)}} \right) \| T \|^{2} 2^{-N/p^{\ast}}.
\end{align*}
So, $\mathbb{P}( \mathcal{A}) <1$ by \eqref{star}. Therefore,  there exists $\theta^{0} = (\theta_{K}^{0})_{K \in \mathcal{B}}$ such that 
$\vert Z(\theta^{0}) \vert < \eta$, $\vert Y_{(m,J)} (\theta^{0}) \vert < \eta$ and 
$\vert W_{(m,J)} (\theta^{0}) \vert < \eta_{(m,J)} $, for every $(m,J) \prec (n,I).$
\end{proof}
\end{lem}

This is the main theorem of this section where we reduce an arbitrary operator to a diagonal operator.
\begin{thm} \label{thm1} Let  $T \colon R_{\omega}^{p,0}  \to R_{\omega}^{p,0}$ be an operator. For every $\epsilon >0$ there exist a diagonal operator  $R \colon R_{\omega}^{p,0}  \to R_{\omega}^{p,0}$  with diagonal elements in  $\mathcal{A}(\mathcal{D}(T))$ such that $T$ is an  orthogonal factor of $R$ with error $\epsilon$.
\begin{proof} Let an arbitrary $\epsilon >0$. Choose an increasing sequence of positive integers $k(n)$ such that 
\begin{align*}
k(n) > p^{\ast} \left( 12n + 13  + 2 \log \left(   \frac{\| T \| } {\epsilon}  \right)  \right) 
\end{align*}
and put $N(n) = k(n) + n$. 

We will construct by induction on the linear order of  $\mathcal{D}_{\omega}$  a  distributional copy $(b^{n}_{I})_{(n,I) \in \mathcal{D}_{\omega}}$ of  $(h^{n}_{I})_{(n,I) \in \mathcal{D}_{\omega}}$ consisting of blocks of the basis $(h^{n}_{I})_{(n,I) \in \mathcal{D}_{\omega}}$    with
\begin{align}
\label{tomato} 
b_{I}^{n} =  b_{\mathcal{B}_{I}^{n}, N(n)} =   \sum _{K \in \mathcal{B}_{I}^{n}} \theta_{K}^{(n,I)} h_{K}^{N(n)},
\end{align}
where  $\mathcal{B}_{I}^{n} \subset \mathcal{D}^{k(n) + i}$ with $I \in \mathcal{D}^{i}$ and $(\theta_{K}^{(n,I)})_{K \in \mathcal{B}_{I}^{n}} \in \lbrace -1, 1 \rbrace ^{\mathcal{B}_{I}^{n}}$. The blocks $(b^{n}_{I})_{(n,I) \in \mathcal{D}_{\omega}}$ will be chosen to satisfy certain orthogonality properties, which will be given explicitly in the inductive construction (see \eqref{first ineq}, \eqref{second ineq} \eqref{third ineq}). The block sequence $(b^{n}_{I})_{(n,I) \in \mathcal{D}_{\omega}}$  will define the operators $j$ and $E$ as in Notation \ref{ntt 2}, and we will show that (see \eqref{this is the end})
\begin{align}
\label{basis small estimates} 
\| (j^{-1} E T j - R)(h_{I}^{n}) \| < \frac{\epsilon}{2^{2n+1 + n/p}},
\end{align}
where $R \colon R_{\omega}^{p,0}  \to  R_{\omega}^{p,0}$ is a diagonal operator defined by 
\begin{align*}
R(h_{I}^{n}) = \left( \vert I \vert ^{-1} \sum_{K \in \mathcal{B}_{I}^{n}} \langle h_{K}^{N(n)}, Th_{K}^{N(n)} \rangle \right) h_{I}^{n} 
\end{align*} 
with diagonal elements in  $\mathcal{A}(\mathcal{D}(T))$. Given \eqref{basis small estimates}, $\| j^{-1} E T j - R \| < \epsilon$ because if $f \in R_{\omega}^{p,0}$ with 
\begin{align*}
f=\sum_{n=1}^{\infty} \sum_{I \in \mathcal{D}_{n-1}} \vert I \vert^{-1} \langle h_{I}^{n}, f \rangle h_{I}^{n},
\end{align*}
then  
\begin{align*}
&\| (j^{-1} E T j - R)(f) \| =  \|  \sum_{n=1}^{\infty} \sum_{I \in \mathcal{D}_{n-1}} \vert I \vert^{-1} \langle h_{I}^{n}, f \rangle  (j^{-1} E T j - R) ( h_{I}^{n}) \|  \\
&\leq \sum_{n=1}^{\infty} \sum_{I \in \mathcal{D}_{n-1}} 2^{n/p} \frac{\epsilon}{2^{n/p + 2n + 1}} \| f \| \leq \epsilon \| f \|.
\end{align*}

Firstly, we start with $n=1$ and $I=[0,1] \in \mathcal{D}^{0}$. 
Take $\mathcal{B}_{I}^{1}  = \mathcal{D}^{k(1)}$. Let
\begin{align*}
b_{I}^{1} (\theta) =   \sum _{K \in \mathcal{B}_{I}^{1}} \theta_{K}^{(1,I)} h_{K}^{N(1)}. 
\end{align*}
and define the random variable 
\begin{align*}
Z_{I}^{1}(\theta) = Z_{\mathcal{B}_{I}^{1}, N(1),T}(\theta) = \langle b_{I}^{1} (\theta), Tb_{I}^{1} (\theta) \rangle - \sum_{K \in \mathcal{B}_{I}^{1} } \langle h_{K}^{N(1)}, Th_{K}^{N(1)}  \rangle,  
\end{align*}
where $\theta = ( \theta^{(1,I)}_{K})_{K \in \mathcal{B}_{I}^{1} } \in \lbrace -1, 1 \rbrace ^{\mathcal{B}_{I}^{1}}$.

By the choice of $k(1)$, $N(1)$ satisfies \eqref{star}, and therefore we may apply Lemma \ref{lem 4}, only to the random variable $Z_I^1$, to choose signs $\theta_{0} \in \lbrace -1, 1 \rbrace ^{\mathcal{B}_{I}^{1}}$ such that 
\begin{align*}
\vert Z_{I}^{1} (\theta_{0}) \vert < \frac{\epsilon}{32^{7}}. 
\end{align*}
We have completed the basis inductive step.

Now, let $(n,I)\in\mathcal D_\omega$ be fixed and assume that we have constructed the blocks $( b_{J}^{m})$ for every $(m,J) \prec (n,I)$ as in \eqref{tomato}. We will only treat the case $I\neq [0,1]$, as the case $I = [0,1]$ is slightly simpler.

We will define a random block
\begin{align*}
b_{I}^{n}(\theta) =  \sum _{K \in \mathcal{B}_{I}^{n}} \theta_{K}^{(n,I)} h_{K}^{N(n)},
\end{align*} 
with $\mathcal{B}_{I}^{n} \subset \mathcal{D}^{k(n) + i}$, to be chosen later, where $i$ is such that $I \in \mathcal{D}^{i}$, such that for the 
random variables
\begin{align*}
Z_{I}^{n} (\theta) & =  Z_{\mathcal{B}_{I}^{n}, N(n),T}(\theta) = \langle b_{I}^{n} (\theta), Tb_{I}^{n} (\theta) \rangle - \sum_{K \in \mathcal{B}_{I}^{n} } \langle h_{K}^{N(n)}, Th_{K}^{N(n)}  \rangle,  
\end{align*}
\begin{align*}
Y^{(n,I)}(m,J)(\theta) & =  Y_{\mathcal{B}_{I}^{n}, N(n), T^{\ast} b_{J}^{m}}(\theta) =  \langle T^{\ast} b_{J}^{m} , b_{I}^{n} (\theta) \rangle  \;\;  \text{for} \;\; (m,J) \prec (n,I) 
\end{align*}
\begin{align*}
W^{(n,I)}(m,J)(\theta) & = W_{\mathcal{B}_{I}^{n}, N(n),Tb_{J}^{m}}(\theta) = \langle b_{I}^{n} (\theta), Tb_{J}^{m}  \rangle  \;\;  \text{for}  \;\; (m,J) \prec (n,I),
\end{align*} 
where $\theta = ( \theta_{K}^{(n,I)})_{ K \in \mathcal{B}_{I}^{n} }  \in \lbrace -1, 1 \rbrace ^{\mathcal{B}_{I}^{n}}$, there exist signs $\theta_{0} \in \lbrace -1, 1 \rbrace ^{\mathcal{B}_{I}^{n}}$ such that  
\begin{align} 
\vert Z_{I}^{n} (\theta_{0} ) \vert & <  \frac{\epsilon}{32^{5n+2}}, \label{first ineq} \\
\vert Y^{(n,I)}(m,J) (\theta_{0} ) \vert & <  \frac{\epsilon}{32^{5n+2}}  \;\; \text{for} \;\;  (m,J) \prec (n,I) \;\; \text{and}  \label{second ineq}\\
 \vert W^{(n,I)}(m,J) (\theta_{0} ) \vert & < \frac{\epsilon}{32^{2n +m +3 +n/q + m/p}}  \;\; \text{for} \;\; (m,J) \prec (n,I). \label{third ineq}
\end{align}
The choice of signs will be given by an appropriate application of Lemma \ref{lem 4}.

We now explicitly define the set $\mathcal{B}_{I}^{n}$. Let $J = \pi(I)$ (such $J$ exists because we assumed $I\neq [0,1]$) and we assume without loss of generality that $I = J^{+}$). So, by the induction hypothesis there exists a block $b_{J}^{n}$ with  $\mathcal{B}_{J}^{n} \subset \mathcal{D}^{k(n) + j}$ and $J \in \mathcal{D}^{j}$. Set 
\begin{align*}
A = \lbrace t \in [0,1] : b_{J}^{n}(t) = 1 \rbrace
\end{align*}
and take  
\begin{align*}
\mathcal{B}_{I}^{n} =  \lbrace K \in \mathcal{D}^{k(n) + i} : \mathrm{supp}(h^{N(n)}_K) \subset A \rbrace. 
\end{align*}
Note that $\mathcal{B}_{I}^{n}$ is the subset of $\mathcal{D}^{k(n)+i}$ such that
\[\cup_{K\in\mathcal{B}_I^n}\mathrm{supp}(h^{N(n)}_K) = A.\]

We define the random block  
\begin{align*}
b_{I}^{n}(\theta) =  \sum _{K \in \mathcal{B}_{I}^{n}} \theta_{K}^{(n,I)} h_{K}^{N(n)},
\end{align*}
where $\theta = ( \theta_{K}^{(n,I)})_{ K \in \mathcal{B}_{I}^{n} }  \in \lbrace -1, 1 \rbrace ^{\mathcal{B}_{I}^{n}}$. We will show that the assumptions of  Lemma \ref{lem 4} are satisfied, and this will yield the existence of sings  $\theta_{0} = (\theta_{K}^{0,(n,I)})  \in \lbrace -1, 1 \rbrace ^{\mathcal{B}_{I}^{n}}$ satisfying inequalities \eqref{first ineq}, \eqref{second ineq} and \eqref{third ineq}.  By the choice of $k(n)$, $N(n)$ satisfies \eqref{star}, and therefore we may apply Lemma \ref{lem 4},  where now the $f_{(m,J)}$ and $x_{(m,J)}$ of  Lemma \ref{lem 4} are $f_{(m,J)} = b_{J}^{m}$ and $x_{(m,J)} = b_{J}^{m}$. By Lemma \ref{lem 4}, the desired signs $\theta_{0} = (\theta_{K}^{0,(n,I)}) $ exist and we put
\begin{align*}
b_{I}^{n}  = b_{I}^{n}(\theta_{0})=   \sum _{K \in \mathcal{B}_{I}^{n}} \theta_{K}^{0,(n,I)} h_{K}^{N(n)}.
\end{align*}
and, for brevity, we denote
\begin{align*}
Z_{I}^{n}  &= Z_{I}^{n} (\theta_{0} )\\
Y^{(n,I)}(m,J) &= Y^{(n,I)}(m,J) (\theta_{0} )\\
W^{(n,I)}(m,J)&= W^{(n,I)}(m,J) (\theta_{0} ).\\
\end{align*}
We have constructed a distributional copy $(b^{n}_{I})_{(n,I) \in \mathcal{D}_{\omega}}$ of  $(h^{n}_{I})_{(n,I) \in \mathcal{D}_{\omega}}$ that satisfies \eqref{first ineq}, \eqref{second ineq} and \eqref{third ineq}.

It remains to show \eqref{basis small estimates}, so fix $(n,I) \in \mathcal{D}_{_{\omega}}$. We have that 
\begin{align*}
&\| (j^{-1} E T j - R)(h_{I}^{n}) \| = \\
& = \|  \sum_{m=1}^{\infty} \sum_{J \in \mathcal{D}_{m-1}}\vert J \vert^{-1} \langle  b_{J}^{m}, Tb_{I}^{n}  \rangle h_{J}^{m} -   \left( \vert I \vert ^{-1} \sum_{K \in \mathcal{B}_{I}^{n}} \langle h_{K}^{N(n)}, Th_{K}^{N(n)} \rangle \right) h_{I}^{n} \| \\
& = \|  \sum_{ (m,J) \prec (n,I)  }  \vert J \vert^{-1} \underbrace{\langle  b_{J}^{m}, Tb_{I}^{n}  \rangle}_{=Y^{(n,I)}(m,J)} h_{J}^{m}   +   \sum_{ (m,J) \succ (n,I)}  \vert J \vert^{-1} \underbrace{\langle  b_{J}^{m}, Tb_{I}^{n}  \rangle}_{=W^{(m,J)}(n,I)} h_{J}^{m} + \\
&  \vert I \vert ^{-1} \underbrace{\left( \langle b_{I}^{n}, T b_{I}^{n} \rangle - \sum_{K \in \mathcal{B}_{I}^{n}} \langle h_{K}^{N(n)}, Th_{K}^{N(n)} \rangle \right)}_{=Z_{I}^{n}} h_{I}^{n} \|. \\
\end{align*}
Therefore,
\begin{align}
\label{this is the end}&\| (j^{-1} E T j - R)(h_{I}^{n}) \|  \\
&\leq 2^{n/q} \vert  Z_{I}^{n} \vert  + 2^{2n+1+n/q} \max_{(m,J) \prec (n,I)}  \vert Y^{(n,I)}(m,J) \vert \nonumber \\
&+ \sum_{m \geq n} 2^{m+1+m/q} \max_{J \in \mathcal{D}_{m-1}} \vert W^{(m,J)}(n,I) \vert < \frac{\epsilon}{2^{2n+1 + n/p}},\nonumber
\end{align}
by the inequalities \eqref{first ineq}, \eqref{second ineq} and \eqref{third ineq}.
\end{proof}

\end{thm}

{Using the above theorem, we prove that the basis of $R_{\omega}^{p,0}$ has the factorization property.} 
\begin{thm} \label{large diagonal theorem} Let  $T \colon R_{\omega}^{p,0}  \to  R_{\omega}^{p,0}$ be an operator with large diagonal 
\begin{align*}
\inf_{(n,I) \in \mathcal{D}_{\omega}} \vert \vert I \vert^{-1} \langle h_{I}^{n}, Th_{I}^{n} \rangle \vert \geq \delta >0.
\end{align*}
Then, for every $0< \epsilon <1$, $T$ is a  $\frac{2(p^{\ast} -1)^{4}}{\delta (1- \epsilon)} (\frac{p^{\ast}}{2})^{3/2}$-factor of the identity operator on  $R_{\omega}^{p,0}$.
\begin{proof} Define the multiplication operator $S \colon R_{\omega}^{p,0}  \to  R_{\omega}^{p,0}$ as the linear extension of 
\begin{align*}
 h_{I}^{n} \longmapsto \left( \vert I \vert^{-1} \langle h_{I}^{n}, Th_{I}^{n} \rangle  \right)^{-1} h_{I}^{n}, \;\;\;  (n,I) \in \mathcal{D}_{\omega}.
 \end{align*}
The operator $S$ is bounded with $\| S \| \leq  (p^{\ast} -1)^{2} / \delta$, where the constant $(p^{\ast} -1)^{2}$ is from the Burkholder inequality, and the operator $TS$ has  diagonal elements equal to $1$.

Let $0< \epsilon < 1$. By Theorem \ref{thm1}, we get that there exists a diagonal operator $R$ such that $TS$ is an orthogonal factor of $R$ with error $\epsilon$, that is $\| j^{-1} E TS j - R \| < \epsilon$, for some $j$ and $E$ as in Notation \ref{ntt 2}.  Because the diagonal entries of $TS$ are one, $R = I$. Hence, $j^{-1}ETSj$ is invertible with $\| \left( j^{-1}ETSj \right)^{-1} \| \leq 1 / (1 - \epsilon)$. Now, put 
\begin{align*}
A &= \left( j^{-1}ETSj \right)^{-1} j^{-1} E 
\end{align*}
and
\begin{align*}
B &= S j. 
\end{align*}
We observe that $ATB = I$. Furthermore, 
\begin{align*}
\| A \| \leq \| \left( j^{-1}ETSj \right)^{-1} \|  \| j^{-1} \| \| E \| \leq \frac{2(p^{\ast} -1)^{2}}{1 - \epsilon} (\frac{p^{\ast}}{2})^{3/2}
\end{align*}
and
\begin{align*}
\| B \| \leq  \| S \| \| j \| \leq \frac{ (p^{\ast} -1)^{2}}{ \delta}.
\end{align*}
Thus, $T$ is a  $\frac{2(p^{\ast} -1)^{4} }{\delta (1- \epsilon)} (\frac{p^{\ast}}{2})^{3/2}$-factor of the identity operator on  $R_{\omega}^{p,0}$.
\end{proof}

\end{thm}

\section{Orthogonal Reduction to a scalar multiple of the identity}

In this section, we approximately orthogonally reduce an arbitrary operator on $R_{\omega}^{p,0}$ to a scalar multiple of the identity on $R_{\omega}^{p,0}$ (Theorem \ref{crucial theorem}) using Proposition \ref{transitivity}, Theorem \ref{thm1} and Theorem \ref{close to a scalar of id}. We use a probabilistic lemma from \cite{lechner:motakis:mueller:schlumprecht:2023} with origin in \cite{lechner:motakis:mueller:schlumprecht:2022} (Lemma \ref{lemma of s4}). Finally, using Theorem \ref{crucial theorem} we prove that for an arbitrary operator $T$ on  $R_{\omega}^{p,0}$, we obtain that either $T$ or $I -T$ is a factor of the identity operator on $R_{\omega}^{p,0}$ (Theorem  \ref{primary fact prop}).

{The proof of the following lemma can be found in \cite[Lemma 3.9]{lechner:motakis:mueller:schlumprecht:2023}. We will use it inside the proof of Theorem \ref{memonomeno theorem} to construct, by induction on the linear order of $\mathcal{D}_{m-1}$,  a distributional copy $(b_{I})_{I \in \mathcal{D}_{m-1}}$ of  $(h_{I})_{I \in \mathcal{D}_{m-1}}$ that stabilizes the eigenvalues of a diagonal operator. As in the proof of Theorem \ref{thm1}, in each step of the induction, we use a random block $b_{I}(\theta)$, defined using a family of dyadic intervals $\mathcal{B}_{I}$ that has been fixed in a previous inductive step. Because $T$ is diagonal, the value $\langle b_I, T (b_I) \rangle$ is independent of the outcome $\theta$, but the latter defines two families of dyadic intervals, $\mathcal B_{I^{+}}(\theta)$ and $\mathcal B_{I^{-}}(\theta)$, to be used in a future inductive step. The role of the Lemma is to establish conditions under which, with positive probability, there is a realization of $b_I$ such that the associated families $\mathcal B_{I^{+}}$ and $\mathcal B_{I^{-}}$, define future random blocks  $b_{I^{+}}(\theta),b_{I^{-}}(\theta)$ with values $\langle b_{I^{+}}, T (b_{I^{+}}) \rangle $ and $\langle b_{I^{-}}, T (b_{I^{-}}) \rangle$ approximately equal.}

\begin{lem} \label{lemma of s4} Let $T \colon L_{p}^{0}  \to  L_{p}^{0}$  be a diagonal operator with  diagonal elements $d_{I} = \vert I \vert^{-1} \langle h_{I}, Th_{I} \rangle$, where $(h_{I})_{I \in \mathcal{D}}$ is  the standard Haar system. Let $\mathcal{B} \subset \mathcal{D}^{m}$ and $m < k$. For signs $\theta = (\theta_{J})_{J \in \mathcal{B}}$ define the  vector 
\begin{align*}
b^{\theta} = \sum_{J \in \mathcal{B}} \theta_{J} h_{J}
\end{align*}   
and denote
\begin{align*}
\Gamma^{+}(\theta) = \lbrace t \in [0,1] : b^{\theta}(t) = 1 \rbrace 
\end{align*}
and 
\begin{align*}
\Gamma^{-}(\theta) = \lbrace t \in [0,1] : b^{\theta}(t) = -1 \rbrace .
\end{align*}
Moreover, denote 
\begin{align*}
\mathcal{B}^{+}(\theta) = \lbrace  J \in \mathcal{D}^{k} : J \subset  \Gamma^{+}(\theta) \rbrace ,
\end{align*}
and 
\begin{align*}
\mathcal{B}^{-}(\theta) = \lbrace  J \in \mathcal{D}^{k} : J \subset  \Gamma^{-}(\theta) \rbrace
\end{align*}
and define the random variables 
\begin{align*}
\lambda^{+} \colon \lbrace -1,1 \rbrace^{\mathcal{B}}  \to  \mathbb{R} \;\;\; \text{and} \;\;\;
\lambda^{-} \colon \lbrace -1,1 \rbrace^{\mathcal{B}}  \to  \mathbb{R} 
\end{align*}
with
\begin{align*}
\lambda^{+}(\theta) = \frac{1}{\text{card}(\mathcal{B}^{+}(\theta))} \sum_{J \in \mathcal{B}^{+}(\theta)} d_{J}
\end{align*}
and
\begin{align*}
\lambda^{-}(\theta) = \frac{1}{\text{card} (\mathcal{B}^{-}(\theta))} \sum_{J \in \mathcal{B}^{-}(\theta)} d_{J}.
\end{align*}

Then 
\begin{align*}
\mathbb{E}(\lambda^{+}) = \mathbb{E}(\lambda^{-}) = \frac{1}{\text{card}( \lbrace J \in \mathcal{D}^{k} : J \subset \cup \mathcal{B}  \rbrace  )} \sum_{ \lbrace J \in \mathcal{D}^{k} : J \subset \cup \mathcal{B}  \rbrace} d_{J},
\end{align*}
\begin{align*}
\mathrm{Var}(\lambda^{+}) \leq \frac{2^{-m}}{\vert \cup \mathcal{B} \vert} \| T \|^{2} 
\end{align*}
and
\begin{align*}
\mathrm{Var}(\lambda^{-}) \leq \frac{2^{-m}}{\vert \cup \mathcal{B} \vert} \| T \|^{2} .
\end{align*}
\end{lem}

\begin{ntt} Let $n \in \mathbb{N}$. For the rest of the paper it is convenient to denote
\begin{align*}
L_{p}^{n} = \langle \lbrace h_{I}^{n} : I \in \mathcal{D}_{n-1} \rbrace \cup \lbrace \chi_{[0,1]} \rbrace \rangle,
\end{align*}
where $ (h_{I}^{n})_{I \in \mathcal{D}_{n-1}}$ is the fixed distributional copy of the standard Haar system $ (h_{I})_{I \in \mathcal{D}_{n-1}}$  used to define $R_\omega^p$. That is, $L_p^n$ is a distributional copy of the space previously denoted by the same symbol. Moreover, we will  denote
\begin{align*}
L_{p}^{n,0} = \lbrace f \in L_{p}^{n} : \mathbb{E}(f) = 0 \rbrace =   \langle \lbrace h_{I}^{n} : I \in \mathcal{D}_{n-1} \rbrace \rangle.
\end{align*}

\end{ntt}

The following notation will be used inside the proofs of Theorems \ref{memonomeno theorem} and \ref{close to a scalar of id}. This is a finite-dimensional version of Notation \ref{ntt 2}.

\begin{ntt} \label{notation for memonomeno theorem } Let $m \leq M$ and $(b_{I})_{I \in \mathcal{D}_{m-1}}$ be a distributional copy of $(h^{m}_{I})_{I \in \mathcal{D}_{m-1}}$ consisting of blocks in the basis $(h^{M}_{I})_{I \in \mathcal{D}_{M-1}}$  of  $L_{p}^{M,0}$. Define the distributional isomorphism $j \colon L_{p}^{m,0}  \to  L_{p}^{M,0}$ with 
\begin{align*}
j(h^{m}_{I}) =  b_{I}, 
\end{align*}
where $I \in \mathcal{D}_{m-1}$ and the projection $E \colon L_{p}^{M,0}  \to  L_{p}^{M,0} $ onto $[(b_{I})]$ defined by
\begin{align*}
E(f) = \sum_{I \in \mathcal{D}_{m-1}} \vert I \vert^{-1} \langle b_{I}, f \rangle b_{I}.
\end{align*}
\end{ntt}

\begin{rem} The projection $E$ is of norm-1 because it is equal to the conditional expectation operator on  $L_{p}^{M,0}$ with respect to $\sigma (b_{I} : I \in \mathcal{D}_{m-1})$.
\end{rem}

{The following finite dimensional reduction to a scalar operator is at the heart of Theorem \ref{close to a scalar of id}. The proof of the latter then relatively easily follows with a compactness argument. To prove this finite dimensional reduction we use a combinatorial argument, based on the pigeonhole principle and a repeated application of the earlier probabilistic future-level stabilization argument (Lemma \ref{lemma of s4}) to stabilize the eigenvalues of a given diagonal operator.}

\begin{thm} \label{memonomeno theorem} Let $m \in \mathbb{N}$, $\Gamma >0$ and $\epsilon >0$. If 
\begin{equation} \label{star gia memonomeno theorem}
M > m+6 +m\frac{4 \Gamma}{\epsilon} + 3 \log(m) + 2 \log(\frac{\Gamma}{\epsilon}) 
\end{equation}
then for every diagonal operator $T \colon L_{p}^{M,0}  \to  L_{p}^{M,0} $ with $\| T \| \leq \Gamma$ and  diagonal elements $d_{I} = \vert I \vert^{-1} \langle h_{I}^{M}, Th_{I}^{M} \rangle$, where $I \in \mathcal{D}_{M-1}$, there exists $\lambda_{0} \in \mathcal{A}(\mathcal{D}(T))$ such that $T$ is an orthogonal factor of $\lambda_{0} I  \colon L_{p}^{m,0}  \to  L_{p}^{m,0} $ with error $(p^{\ast} -1)  \epsilon $, where $I$ is the identity operator on $L_{p}^{m,0}$.
\begin{proof} Let 
\begin{align*}
x = [m+6 + 3 \log(m) + 2 \log(\frac{\Gamma}{\epsilon})] +1.
\end{align*}
For $k \in \lbrace x, x+1,...,M-1 \rbrace$, define 
\begin{align*}
\lambda_{k} = \frac{1}{\text{card}(\mathcal{D}^{k})} \sum_{J \in \mathcal{D}^{k}} d_{J} 
\end{align*}
and observe that $\vert \lambda_{k} \vert \leq \| T \|$, for every  $k \in \lbrace x, x+1,...,M-1 \rbrace$. By \eqref{star gia memonomeno theorem} and the pigeonhole principle there exist 
\begin{align*}
x \leq n_{1}<n_{2} <...<n_{m-1} \leq M -1
\end{align*} 
such that 
\begin{align*}
\vert \lambda_{n_{k}}  -  \lambda_{n_{l}} \vert  < \frac{\epsilon}{2}, 
\end{align*}
for every $1 \leq k,l \leq m-1$. 

We construct by induction on the linear order of $\mathcal{D}_{m-1}$ families $(\mathcal{B}_{I})_{I \in \mathcal{D}_{m-1}}$ and signs $((\theta_{J})_{J \in \mathcal{B}_{I}})_{I \in  \mathcal{D}_{m-1}}$ with $\mathcal{B}_{[0,1]} = \mathcal{D}^{n_{1}}$ and if 
\begin{align*}
b_{I} = \sum_{J \in \mathcal{B}_{I}} \theta_{J} h_{J}^{M}, 
\end{align*}
where $I \in \mathcal{D}^{k}$, then 
\begin{align*}
\mathcal{B}_{I^{+}} = \lbrace J \in   \mathcal{D}^{n_{k+1}} : J \subset [b_{I} = 1] \rbrace
\end{align*}
and
\begin{align*}
\mathcal{B}_{I^{-}} = \lbrace J \in   \mathcal{D}^{n_{k+1}} : J \subset [b_{I} = -1] \rbrace.
\end{align*}
From these, we get a distributional copy $(b_{I})_{I \in \mathcal{D}_{m-1}}$ of $(h^{m}_{I})_{I \in \mathcal{D}_{m-1}}$ that defines the operators $j$ and $E$ in Notation \ref{notation for memonomeno theorem } and $(\lambda_{I})_{I \in \mathcal{D}_{m-1}} \subset \mathcal{A}(\mathcal{D}(T))$ with $\vert \lambda_{I} - \lambda_{J} \vert < \epsilon$, for every $I,J \in \mathcal{D}_{m-1}$, such that $\| j^{-1} E T j - \lambda_{0} I \| < (p^{\ast} -1)  \epsilon $, where $\lambda_{0} \in  \mathcal{A}(\mathcal{D}(T))$.

We begin the inductive construction. Although the first inductive step is similar to the general one, we include for completeness. Let $\mathcal{B}_{[0,1]} = \mathcal{D}^{n_{1}}$ and for signs $\theta = (\theta_{J})_{J \in \mathcal{B}_{[0,1]}}$ define the vector 
\begin{align*}
b_{[0,1]}^{\theta} = \sum_{J\in \mathcal{B}_{[0,1]}} \theta_{J} h_{J}^{M}
\end{align*}
and denote 
\begin{align*}
\Gamma^{+}_{[0,1]}(\theta) = \lbrace t \in [0,1] : b^{\theta}_{[0,1]}(t) = 1 \rbrace 
\end{align*}
and
\begin{align*}
\Gamma^{-}_{[0,1]}(\theta) = \lbrace t \in [0,1] : b^{\theta}_{[0,1]}(t) = -1 \rbrace .
\end{align*}
Moreover, denote
\begin{align*}
\mathcal{B}_{[0,\frac{1}{2}]}^{n_{k}} (\theta) = \lbrace  J \in \mathcal{D}^{n_{k}} : J \subset   \Gamma^{+}_{[0,1]}(\theta)  \rbrace ,
\end{align*}
\begin{align*}
\mathcal{B}_{[\frac{1}{2},1]}^{n_{k}} (\theta) = \lbrace  J \in \mathcal{D}^{n_{k}} : J \subset   \Gamma^{-}_{[0,1]}(\theta)  \rbrace , 
\end{align*}
for every $2 \leq k \leq m-1$, and define the  random variables 
\begin{align*}
\lambda_{[0,\frac{1}{2}]}^{n_{k}} \colon \lbrace -1,1 \rbrace^{\mathcal{B}}  \to  \mathbb{R} \;\;\; \text{and} \;\;\; \lambda_{[\frac{1}{2},1]}^{n_{k}} \colon \lbrace -1,1 \rbrace^{\mathcal{B}}  \to  \mathbb{R} 
\end{align*}
with
\begin{align*}
\lambda_{[0,\frac{1}{2}]}^{n_{k}}(\theta) = \frac{1}{\text{card} ( \mathcal{B}_{[0,\frac{1}{2}]}^{n_{k}} (\theta) ) } \sum_{J \in  \mathcal{B}_{[0,\frac{1}{2}]}^{n_{k}} (\theta) } d_{J}
\end{align*}
and
\begin{align*}
\lambda_{[\frac{1}{2},1]}^{n_{k}}(\theta) = \frac{1}{\text{card} ( \mathcal{B}_{[\frac{1}{2},1]}^{n_{k}} (\theta) ) } \sum_{J \in  \mathcal{B}_{[\frac{1}{2},1]}^{n_{k}} (\theta) } d_{J},
\end{align*}
for every $2 \leq k \leq m-1$. We are looking for signs $\theta = (\theta_{J})_{J \in \mathcal{B}_{[0,1]}}$  such that 
\begin{align*}
\vert \lambda_{[0,\frac{1}{2}]}^{n_{k}}(\theta) - \lambda_{n_{k}} \vert < \frac{\epsilon}{4m}
\end{align*} 
and 
\begin{align*}
\vert \lambda_{[\frac{1}{2},1]}^{n_{k}}(\theta) - \lambda_{n_{k}} \vert < \frac{\epsilon}{4m}
\end{align*} 
for every $2 \leq k \leq m-1$. But, by Lemma \ref{lemma of s4}, we have that 
\begin{align*}
\mathbb{E}(\lambda_{[0,\frac{1}{2}]}^{n_{k}}) = \mathbb{E}(\lambda_{[\frac{1}{2},1]}^{n_{k}}) = \lambda_{n_{k}},
\end{align*}
\begin{align*}
\mathrm{Var}(\lambda_{[0,\frac{1}{2}]}^{n_{k}}) \leq \frac{1}{2^{n_{1}}} \Gamma^{2}
\end{align*}
and
\begin{align*}
\mathrm{Var}(\lambda_{[\frac{1}{2},1]}^{n_{k}}) \leq \frac{1}{2^{n_{1}}} \Gamma^{2},
\end{align*}
for every $2 \leq k \leq m-1$. Now, define the events, 
\begin{align*}
\mathcal{Z}_{[0,\frac{1}{2}]}^{n_{k}} & = \lbrace \theta : \vert \lambda_{[0,\frac{1}{2}]}^{n_{k}}(\theta) - \lambda_{n_{k}} \vert > \frac{\epsilon}{4m} \rbrace \;\; \text{and}\\
\mathcal{Z}_{[\frac{1}{2},1]}^{n_{k}} & = \lbrace \theta : \vert \lambda_{[\frac{1}{2},1]}^{n_{k}}(\theta) - \lambda_{n_{k}} \vert > \frac{\epsilon}{4m} \rbrace, 
\end{align*}
for every $2 \leq k \leq m-1$, and by the  Chebyshev inequality we get that 
\begin{align*}
\mathbb{P} (\mathcal{Z}_{[0,\frac{1}{2}]}^{n_{k}}) \leq \frac{16m^{2}}{\epsilon^{2}} \mathrm{Var}(\lambda_{[0,\frac{1}{2}]}^{n_{k}}) \leq 16m^{2} \frac{\Gamma^{2}}{\epsilon^{2}} \frac{1}{2^{n_{1}}}
\end{align*}
and
\begin{align*}
\mathbb{P} (\mathcal{Z}_{[\frac{1}{2},1]}^{n_{k}}) \leq \frac{16m^{2}}{\epsilon^{2}} \mathrm{Var}(\lambda_{[\frac{1}{2},1]}^{n_{k}}) \leq 16m^{2} \frac{\Gamma^{2}}{\epsilon^{2}} \frac{1}{2^{n_{1}}},
\end{align*} 
for every $2 \leq k \leq m-1$. So, 
\begin{align*}
&\mathbb{P} \left(   \bigcup_{2 \leq k \leq m}   \mathcal{Z}_{[0,\frac{1}{2}]}^{n_{k}} \cup   \mathcal{Z}_{[\frac{1}{2},1]}^{n_{k}}  \right) \leq \sum_{2 \leq k \leq m} \mathbb{ P} ( \mathcal{Z}_{[0,\frac{1}{2}]}^{n_{k}}) + \mathbb{P} (\mathcal{Z}_{[\frac{1}{2},1]}^{n_{k}})\\
&\leq \sum_{2 \leq k \leq m } 32m^{2} \frac{\Gamma^{2}}{\epsilon^{2}} \frac{1}{2^{n_{1}}} \leq  2^{5} m^{3}  \frac{\Gamma^{2}}{\epsilon^{2}} \frac{1}{2^{n_{1}}} <1,
\end{align*}
because $x \leq n_{1}$. Hence,  there exist signs $\theta = (\theta_{J})_{J \in \mathcal{B}_{[0,1]}}$  such that 
\begin{align*}
\vert \lambda_{[0,\frac{1}{2}]}^{n_{k}}(\theta) - \lambda_{n_{k}} \vert < \frac{\epsilon}{4m}
\end{align*} 
and 
\begin{align*}
\vert \lambda_{[\frac{1}{2},1]}^{n_{k}}(\theta) - \lambda_{n_{k}} \vert < \frac{\epsilon}{4m}
\end{align*} 
for every $2 \leq k \leq m-1$. 

Now, assume that for $I\in \mathcal{D}^{k-1}$ we have found signs $\theta^{\prime} = (\theta^{\prime}_{J})_{J\in \mathcal{B}_{\pi (I)}}$, where without loss of generality assume that $I = \pi (I)^{+}$, sets 
\begin{align*}\mathcal{B}_{I}^{n_{j}} = \lbrace J \in \mathcal{D}^{n_{j}} : J \subset [b_{\pi (I)}^{\theta^{\prime}} = 1] \rbrace 
\end{align*}
and 
\begin{align*}
\lambda_{I}^{n_{j}} = \frac{1}{\text{card} (\mathcal{B}_{I}^{n_{j}})   } \sum_{J \in \mathcal{B}_{I}^{n_{j}}} d_{J} 
\end{align*} 
such that 
\begin{align*}
\vert \lambda_{I}^{n_{j}}  -  \lambda_{\pi(I)}^{n_{j}} \vert < \frac{\epsilon}{4m} 
\end{align*}
for every $k \leq j \leq m-1$. Now, for $I^{+},I^{-} \in \mathcal{D}^{k}$, we are looking for signs $\theta = (\theta_{J})_{J \in \mathcal{B}_{I}^{n_{k}}}$, such that for the sets 
\begin{align*}
\mathcal{B}_{I^{+}}^{n_{j}} = \lbrace J \in \mathcal{D}^{n_{j}} : J \subset [b_{I}^{\theta} = 1] \rbrace 
\end{align*}
and  
\begin{align*}
\mathcal{B}_{I^{-}}^{n_{j}} = \lbrace J \in \mathcal{D}^{n_{j}} : J \subset [b_{I}^{\theta} = -1] \rbrace, 
\end{align*}
we have  
\begin{align*}
\lambda_{I^{+}}^{n_{j}} = \frac{1}{\text{card} (\mathcal{B}_{I^{+}}^{n_{j}})   } \sum_{J \in \mathcal{B}_{I}^{n_{j}}} d_{J}
\end{align*}
and 
\begin{align*} 
\lambda_{I^{-}}^{n_{j}} = \frac{1}{\text{card} (\mathcal{B}_{I^{-}}^{n_{j}})   } \sum_{J \in \mathcal{B}_{I}^{n_{j}}} d_{J}, 
\end{align*}
such that 
\begin{align*}
\vert \lambda_{I^{+}}^{n_{j}}  -  \lambda_{I}^{n_{j}} \vert < \frac{\epsilon}{4m}
\end{align*}
and
\begin{align*}
\vert \lambda_{I^{-}}^{n_{j}}  -  \lambda_{I}^{n_{j}} \vert < \frac{\epsilon}{4m} 
\end{align*}
for every $k+1 \leq j \leq m-1$. 
Indeed, let  $\theta = (\theta_{J})_{J \in \mathcal{B}_{I}^{n_{k}}}$. For the vector 
\begin{align*}
b_{I}^{\theta} = \sum_{J \in \mathcal{B}_{I}^{n_{k}}} \theta_{J} h_{J}^{M},
\end{align*}
denote 
\begin{align*}
\Gamma^{+}_{I}(\theta) = \lbrace t \in [0,1] : b^{\theta}_{I}(t) = 1 \rbrace 
\end{align*}
and
\begin{align*}
\Gamma^{-}_{I}(\theta) = \lbrace t \in [0,1] : b^{\theta}_{I}(t) = -1 \rbrace .
\end{align*}
Moreover, denote
\begin{align*}
\mathcal{B}_{I^{+}}^{n_{j}} (\theta) = \lbrace  J \in \mathcal{D}^{n_{j}} : J \subset   \Gamma^{+}_{I}(\theta)  \rbrace ,
\end{align*}
\begin{align*}
\mathcal{B}_{I^{-}}^{n_{j}} (\theta) = \lbrace  J \in \mathcal{D}^{n_{j}} : J \subset   \Gamma^{-}_{I}(\theta)  \rbrace ,
\end{align*}  
and define the  random variables 
\begin{align*}
\lambda_{I^{+}}^{n_{j}} \colon \lbrace -1,1 \rbrace^{\mathcal{B}_{I}^{n_{k}}}  \to  \mathbb{R} \;\;\; \text{and} \;\;\; \lambda_{I^{-}}^{n_{j}} \colon \lbrace -1,1 \rbrace^{\mathcal{B}_{I}^{n_{k}}}  \to  \mathbb{R} 
\end{align*}
with 
\begin{align*}
\lambda_{I^{+}}^{n_{j}}(\theta) = \frac{1}{\text{card} (\mathcal{B}_{I^{+}}^{n_{j}} (\theta))   } \sum_{J \in \mathcal{B}_{I^{+}}^{n_{j}}(\theta)} d_{J}
\end{align*}
and 
\begin{align*}
\lambda_{I^{-}}^{n_{j}}(\theta) = \frac{1}{\text{card} (\mathcal{B}_{I^{-}}^{n_{j}} (\theta))   } \sum_{J \in \mathcal{B}_{I^{-}}^{n_{j}}(\theta)} d_{J}
\end{align*}
We are looking for signs $\theta = (\theta_{J})_{J \in \mathcal{B}_{I}^{n_{k}}}$ such that 
\begin{align*}
\vert \lambda_{I^{+}}^{n_{j}}(\theta)  -  \lambda_{I}^{n_{j}} \vert < \frac{\epsilon}{4m}
\end{align*} 
and 
\begin{align*}
\vert \lambda_{I^{-}}^{n_{j}}(\theta)  -  \lambda_{I}^{n_{j}} \vert < \frac{\epsilon}{4m} 
\end{align*}
for every $k+1 \leq j \leq m-1$. Indeed, define the events 
\begin{align*}
\mathcal{Z}_{I^{+}}^{n_{j}} = \lbrace \theta :  \vert \lambda_{I^{+}}^{n_{j}}(\theta)  -  \lambda_{I}^{n_{j}} \vert > \frac{\epsilon}{4m}  \rbrace 
\end{align*}
and
\begin{align*} 
\mathcal{Z}_{I^{-}}^{n_{j}} = \lbrace \theta :  \vert \lambda_{I^{-}}^{n_{j}}(\theta)  -  \lambda_{I}^{n_{j}} \vert > \frac{\epsilon}{4m}  \rbrace ,
\end{align*} 
for every $k+1 \leq j \leq m-1$. By Lemma \ref{lemma of s4}, we get 
\begin{align*}
\mathbb{E}(\lambda_{I^{+}}^{n_{j}}) = \mathbb{E}(\lambda_{I^{-}}^{n_{j}}) =  \lambda_{I}^{n_{j}} ,
\end{align*}
\begin{align*}
\mathrm{Var}(  \lambda_{I^{+}}^{n_{j}}) \leq \frac{2^{k-1}}{2^{n_{k}}}  \Gamma^{2}  
\end{align*}
and
\begin{align*}
\mathrm{Var}(  \lambda_{I^{-}}^{n_{j}}) \leq \frac{2^{k-1}}{2^{n_{k}}}  \Gamma^{2},  
\end{align*}
for every $k+1 \leq j \leq m-1$. But, by the Chebyshev inequality we get that 
\begin{align*}
\mathbb{P}(   \mathcal{Z}_{I^{+}}^{n_{j}}) \leq 16m^{2} \epsilon^{-2}  \mathrm{Var}(  \lambda_{I^{+}}^{n_{j}})
\end{align*}
and  
\begin{align*}
\mathbb{P}(   \mathcal{Z}_{I^{-}}^{n_{j}}) \leq 16m^{2} \epsilon^{-2}  \mathrm{Var}(  \lambda_{I^{-}}^{n_{j}}),
\end{align*}
for every $k+1 \leq j \leq m-1$. Therefore, 
\begin{align*}
\mathbb{P} \left(  \bigcup_{j=k+1}^{m}    \mathcal{Z}_{I^{+}}^{n_{j}} \cup     \mathcal{Z}_{I^{-}}^{n_{j}} \right) \leq \sum_{j=k+1}^{m} \mathbb{P}(\mathcal{Z}_{I^{+}}^{n_{j}}) +\mathbb{P}( \mathcal{Z}_{I^{-}}^{n_{j}}) \\ 
\leq \sum_{j=k+1}^{m} 2^{5}  \frac{m^{2}}{\epsilon^{2}} \frac{2^{k-1}}{2^{n_{k}}} \Gamma^{2}
< 2^{m+5} m^{3} \frac{\Gamma^{2}}{\epsilon^{2}} \frac{1}{2^{n_{k}}} <1,
\end{align*}
because $x \leq n_{1}<n_{2} <...<n_{m-1}$. Hence, there exist $\theta = (\theta_{J})_{J \in \mathcal{B}_{I}^{n_{k}}}$ such that 
\begin{align*}
\vert \lambda_{I^{+}}^{n_{j}}(\theta)  -  \lambda_{I}^{n_{j}} \vert < \frac{\epsilon}{4m}
\end{align*}
and 
\begin{align*}
\vert \lambda_{I^{-}}^{n_{j}}(\theta)  -  \lambda_{I}^{n_{j}} \vert < \frac{\epsilon}{4m}, 
\end{align*}
for every $k+1 \leq j \leq m-1$.

Now, take $\lambda_{0} = \lambda_{n_{1}}$. For $I \in \mathcal{D}^{k}$ and $k+1 \leq j \leq m-1$,  we observe that 
\begin{align*}
\vert \lambda_{I}^{n_{j}} - \lambda_{n_{j}} \vert \leq \vert \sum_{i=0}^{k-1} \lambda_{\pi(I)^{i}}^{n_{j}} -  \lambda_{\pi(I)^{i+1}}^{n_{j}} \vert < \frac{\epsilon}{4}, 
\end{align*}
for every $k+1 \leq j \leq m-1$, where $\lambda_{\pi(I)^{0}}^{n_{j}} = \lambda_{I}^{n_{j}}$ and $\lambda_{\pi(I)^{k}}^{n_{j}} = \lambda_{n_{j}}$. Therefore, 
\begin{align*}
\vert \lambda_{I}^{n_{j}}  - \lambda_{0} \vert \leq \vert \lambda_{I}^{n_{j}}  - \lambda_{n_{j}} \vert  +   \vert \lambda_{n_{j}}  - \lambda_{0} \vert < \frac{\epsilon}{4} + \frac{\epsilon}{2} < \epsilon ,
\end{align*}
for every $k+1 \leq j \leq m-1$. Now, for $I \in \mathcal{D}^{k}$, put $\lambda_{I} = \lambda_{I}^{n_{k+1}}$. We have that $\vert \lambda_{I} -   \lambda_{0}  \vert < \epsilon,$ for every  $I \in \mathcal{D}_{m-1}$. Finally, for an arbitrary $f \in L_{p}^{m,0}$ with $f = \sum_{I \in \mathcal{D}_{m-1}} a_{I} h_{I}^{m}$, by the Paley-Burkholder inequality, we get that 
\begin{align*}
& \|(j^{-1} E T j - \lambda_{0} I)(f) \| = \| \sum_{I \in \mathcal{D}_{m-1}} (\lambda_{I} - \lambda_{0}) a_{I} h_{I}^{m} \| \\
&\leq (p^{\ast} -1) \epsilon \|\sum_{I \in \mathcal{D}_{m-1}} a_{I} h_{I}^{m} \| = (p^{\ast} -1)  \epsilon \|f\|.
\end{align*}
Therefore, $\| j^{-1} E T j - \lambda_{0} I \| < (p^{\ast} -1)  \epsilon$.
\end{proof}

\end{thm}

We now reduce an arbitrary diagonal operator to a scalar multiple of the identity operator using the above finite dimensional reduction.
\begin{thm} \label{close to a scalar of id}  Let  $T \colon R_{\omega}^{p,0}  \to R_{\omega}^{p,0}$ be a diagonal operator with diagonal elements $d_{I}^{n} = \vert I \vert^{-1} \langle h_{I}^{n}, Th_{I}^{n} \rangle$, where $(n,I) \in \mathcal{D}_{\omega}$. For every $\epsilon >0$ there exists $\lambda_{0} \in \mathcal{A}(\mathcal{D}(T))$ such that $T$ is an orthogonal factor of $\lambda_{0}I  \colon R_{\omega}^{p,0}  \to R_{\omega}^{p,0} $ with error $\epsilon$, where $I$ is the identity operator on  $R_{\omega}^{p,0}$.
\begin{proof} For $\frac{\epsilon}{4 (p^{\ast} -1) }>0$, choose an increasing sequence of non-negative integers $N(m)$ with 
\begin{align*}
N(m)> m+6 + m\frac{16 (p^{\ast} -1)   \| T \|}{\epsilon}  + 3 \log(m) +2 \log(\frac{ 4 (p^{\ast} -1)  \| T \|}{\epsilon}).
\end{align*}
Then by Theorem \ref{memonomeno theorem}, there exist $\lambda_{0}^{m} \in  \mathcal{A}(\mathcal{D}(T))$, where 
\begin{align*}
\vert \lambda_{I}^{m} - \lambda_{J}^{m}  \vert < \frac{\epsilon}{4 (p^{\ast} -1) },
\end{align*} 
for every $ I,J \in  \mathcal{D}_{m-1} $, and a distributional copy $(b_{I}^{m})_{I \in \mathcal{D}_{m-1}}$ of $(h^{m}_{I})_{I \in \mathcal{D}_{m-1}}$ consisting of blocks in the basis $(h^{N(m)}_{I})_{I \in \mathcal{D}_{N(m)-1}}$ of $L_{p}^{N(m),0}$ that defines the operators $j_{m}$ and $E_{m}$ in Notation  \ref{notation for memonomeno theorem } such that 
\begin{align*}
\| j_{m}^{-1} E_{m} T j_{m} - \lambda_{0}^{m} I_{m} \| <  \frac{\epsilon}{4}. \end{align*} 
Since the sequence $(\lambda_{0}^{m})$ is bounded there exists a subsequence  $(\lambda_{0}^{m_{k}})$ and $\lambda_{0} \in \mathcal{A}(\mathcal{D}(T))$ such that 
\begin{align*}
\vert \lambda_{0}^{m_{k}} - \lambda_{0}  \vert < \frac{\epsilon}{2 (p^{\ast} -1) }
\end{align*} 
for every $k=1,2...$. Thus, 
\begin{align*}
\vert \lambda_{I}^{m_{k}} - \lambda_{0} \vert < \frac{\epsilon}{(p^{\ast} -1)  },
\end{align*}
for every $k=1,2...$. Now, put 
\begin{align*}
Y = [ (b_{I}^{m_{k}})_{I \in \mathcal{D}_{k-1}} : k=1,2,...].
\end{align*}
Define the distributional isomorphism $j \colon R_{\omega}^{p,0}  \to Y$ with $j(h_{I}^{k}) = b_{I}^{m_{k}}$, where $I \in \mathcal{D}_{k-1}$, and the projection $E \colon R_{\omega}^{p,0}  \to R_{\omega}^{p,0}$ onto $Y$ with 
\begin{align*}
E(f) = \sum_{k=1}^{\infty} \sum_{I \in \mathcal{D}_{k-1}} \vert I \vert^{-1} \langle b_{I}^{m_{k}}, f \rangle   b_{I}^{m_{k}},
\end{align*}
which is bounded by Remark \ref{important for boundedness of E}.  Finally, because $j^{-1} E T j$ is a diagonal operator with entries $(\lambda_{I}^{m_{k}})_{(k,I)\in \mathcal{D}_{\omega}} $ we get that 
\begin{align*}
\| j^{-1} E T j - \lambda_{0}I \| < \epsilon .
\end{align*}
\end{proof}

\end{thm}

In the following theorem, we reduce an arbitrary operator to a scalar multiple of the identity operator.
\begin{thm} \label{crucial theorem} Let $T \colon R_{\omega}^{p,0}  \to R_{\omega}^{p,0}$ be an operator. For every $\epsilon >0$ there exists $\lambda_{0}$ in the convex hull of  $\mathcal{D}(T)$ such that $T$ is an orthogonal factor of $\lambda_{0}I  \colon R_{\omega}^{p,0}  \to R_{\omega}^{p,0} $ with error $\epsilon$, where $I$ is the identity operator on  $R_{\omega}^{p,0}$.
\begin{proof} It follows immediately from Proposition \ref{transitivity}, Theorem \ref{thm1} and Theorem \ref{close to a scalar of id}.
\end{proof}

\end{thm}

Using the above theorem, we prove that $ R_{\omega}^{p,0}$ has the primary factorization property.
\begin{thm} \label{primary fact prop} Let $T \colon R_{\omega}^{p,0}  \to R_{\omega}^{p,0}$ be an operator. For every $0< \epsilon < 1$, either $T$ or $I -T$ is a $ \frac{4}{1 -  \epsilon} (p^{\ast} -1 )^{2} (\frac{p^{\ast}}{2})^{3/2}$-factor of  the identity operator $I$ on  $R_{\omega}^{p,0}$.
\begin{proof} Let $0< \epsilon < 1$. By Theorem \ref{crucial theorem}, there exists $\lambda_{0}$ in the convex hull of  $\mathcal{D}(T)$ such that $T$ is an orthogonal factor of $\lambda_{0}I  \colon R_{\omega}^{p,0}  \to R_{\omega}^{p,0} $ with error $\epsilon /2$, that is there exists a distributional embedding $j \colon R_{\omega}^{p,0}  \to R_{\omega}^{p,0}$ such that 
\begin{align*}
\| j^{-1} P_{j(R_{\omega}^{p,0})}Tj - \lambda_{0} I \| < \epsilon / 2.
\end{align*}
If $\vert \lambda_{0}  \vert \geq 1/2$, then 
\begin{align*}
\| \lambda_{0}^{-1} j^{-1} P_{j(R_{\omega}^{p,0})} T j - I \|  <  \epsilon .
\end{align*}
So, the operator $A = \lambda_{0}^{-1} j^{-1} P_{j(R_{\omega}^{p,0})} T j$ is invertible with $\| A^{-1} \| \leq 1/(1 - \epsilon)$. Moreover, if $B = A^{-1} \lambda_{0}^{-1} j^{-1}  P_{j(R_{\omega}^{p,0})}$, then $BTj = I$ and $\| B \| \| j \| \leq 2/(1- \epsilon) \| P_{j(R_{\omega}^{p,0})} \| \leq 4/(1- \epsilon) (p^{\ast} -1)^{2} (\frac{p^{\ast}}{2})^{3/2}$. Now, if $\vert \lambda_{0}  \vert < 1/2$, then, because 
\begin{align*}
\| j^{-1} P_{j(R_{\omega}^{p,0})} ( I - T) j - (1- \lambda_{0}) I \| < \epsilon ,
\end{align*}
we achieve the same conclusion for $I-T$ instead of $T$.
\end{proof}

\end{thm}

\subsection*{Acknowledgement} We are thankful to the anonymous referees for providing valuable suggestions to improve the exposition of our paper.


%
%
\bibliographystyle{plain}
\bibliography{bibliography}

\end{document}